\RequirePackage{fix-cm}
\documentclass[smallextended]{svjour3}       % onecolumn (second format)
\smartqed  % flush right qed marks, e.g. at end of proof
\usepackage{graphicx}
\usepackage{amsmath}
\usepackage{amssymb}
\usepackage{latexsym}
\usepackage{mathrsfs}
\usepackage{colortbl}
\usepackage{amscd}
\usepackage{tikz-cd}
\usepackage[all]{xy}
\usepackage{comment}

 \journalname{Japan Journal of Industrial and Applied Mathematics}

\def\diam{\mathop{\mathrm{diam}}}
\def\diag{\mathop{\mathrm{diag}}}
\def\>{\textgreater}
\def\<{\textless}
\def\Span{\mathop{\mathrm{span}}}

\def\conv{\mathop{\mathrm{conv}}}

\def\Int{\mathop{\rm{int}}}
\def\Abs{\mathop{\rm{abs}}}

\spnewtheorem{thr}{Theorem}{\bf}{\it}
\spnewtheorem{coro}{Corollary}{\bf}{\it}
\spnewtheorem{defi}{Definition}{\bf}{\it}
\spnewtheorem{lem}{Lemma}{\bf}{\it}
\spnewtheorem{prop}{Proposition}{\bf}{\it}
\spnewtheorem{ass}{Assumption}{\bf}{\it}
\spnewtheorem{cond}{Condition}{\bf}{\it}
\spnewtheorem{Rem}{Remark}{\it}{\it}
\spnewtheorem{Ex}{Example}{\it}{\it}
\spnewtheorem{Note}{Note}{\it}{\it}
\spnewtheorem*{ex*}{Example:}{\it}{\it}
\spnewtheorem*{pf*}{Proof}{\bf}{\rm}
\spnewtheorem*{rem*}{Remark:}{\it}{\it}
\spnewtheorem*{note*}{Note:}{\it}{\it}

\spnewtheorem*{lem1*}{Lemma 1}{\bf}{\rm}
\spnewtheorem*{lem3*}{Lemma 3}{\bf}{\rm}

\spnewtheorem*{th21*}{Theorem A}{\bf}{\it}
\spnewtheorem*{thrBA*}{Theorem B}{\bf}{\it}
\spnewtheorem*{thrF*}{Theorem C}{\bf}{\it}
\spnewtheorem*{thrin*}{Theorem D}{\bf}{\it}
\spnewtheorem*{thrlp*}{Theorem E}{\bf}{\it}
\spnewtheorem*{thrso*}{Theorem}{\bf}{\it}

\spnewtheorem*{thA0*}{Theorem A0}{\bf}{\it}

\spnewtheorem*{th2*}{Theorem 2}{\bf}{\it}

\spnewtheorem*{coex1*}{Counterexample 1}{\it}{\it}
\spnewtheorem*{coex2*}{Counterexample 2}{\it}{\it}

\allowdisplaybreaks[3]

\newcounter{sone}
\newcounter{stwo}
\newcounter{sthree}
\newcounter{sfour}
\newcounter{sfive}
\newcounter{ssix}
\newcounter{lone}
\newcounter{ltwo}
\newcounter{lthree}
\newcounter{lfour}
\newcounter{lfive}
\newcounter{lsix}
\begin{document}

%\title{General theory of interpolation error estimates on anisotropic meshes, corrigendum}
\title{Anisotropic interpolation error estimates using a new geometric parameter}
%\title{General theory of interpolation error estimates on anisotropic meshes, part \Roman{ltwo}}%\thanks{Grants or other notes
%about the article that should go on the front page should be
%placed here. General acknowledgments should be placed at the end of the article.}}
%\title{Remarks on interpolation error analysis and inverse inequalities on anisotropic meshes}
%\title{Remarks on anisotropic interpolation error analysis}
%\title{Classification of interpolation error estimates on anisotropic meshes}
%\subtitle{Do you have a subtitle?\\ If so, write it here}

\titlerunning{Interpolation error analysis on anisotropic meshes}        % if too long for running head

\author{Hiroki Ishizaka \and Kenta Kobayashi \and Takuya Tsuchiya %etc.
}

%\authorrunning{Short form of author list} % if too long for running head

\institute{Hiroki Ishizaka \at
              Graduate School of Science and Engineering, Ehime University, Matsuyama, Japan \\
              \email{h.ishizaka005@gmail.com}           %  \\
%             \emph{Present address:} of F. Author  %  if needed
           \and
           Kenta Kobayashi\at
              Graduate School of Business Administration, Hitotsubashi University, Kunitachi, Japan \\
            \email{kenta.k@r.hit-u.ac.jp}
            \and
            Takuya Tsuchiya \at 
             Graduate School of Science and Engineering, Ehime University, Matsuyama, Japan \\
              \email{tsuchiya@math.sci.ehime-u.ac.jp}  
}

\date{Received: date / Accepted: date}
% The correct dates will be entered by the editor

\maketitle

\begin{abstract}
%We present precise anisotropic interpolation error estimates for smooth functions using a new geometric parameter and inverse inequalities on anisotropic meshes. In our theory, the error of interpolation is bound in terms of the diameter of a simplex and the geometric parameter. Imposing additional assumptions makes it possible to obtain anisotropic error estimates. This paper also includes corrections to an error in "General theory of interpolation error estimates on anisotropic meshes" (Japan Journal of Industrial and Applied Mathematics, 38 (2021) 163-191), in which Theorem 2 is incorrect.  

%We present precise anisotropic interpolation error estimates for smooth functions using a new geometric parameter and derive inverse inequalities on anisotropic meshes. In our theory, the interpolation error is bounded by a function of the diameter of a simplex and the geometric parameter. Imposing additional assumptions makes it possible to obtain anisotropic error estimates. This paper also includes corrections to an error in Theorem 2 of our previous paper, ``General theory of interpolation error estimates on anisotropic meshes" (Japan Journal of Industrial and Applied Mathematics, 38 (2021) 163-191). 

We present precise anisotropic interpolation error estimates for smooth functions using a new geometric parameter and derive inverse inequalities on anisotropic meshes. In our theory, the interpolation error is bounded in terms of the diameter of a simplex and the geometric parameter. Imposing additional assumptions makes it possible to obtain anisotropic error estimates. This paper also includes corrections to an error in Theorem 2 of our previous paper, ``General theory of interpolation error estimates on anisotropic meshes" (Japan Journal of Industrial and Applied Mathematics, 38 (2021) 163-191). 

\keywords{Finite element \and Interpolation error estimates \and Anisotropic meshes}
% \PACS{PACS code1 \and PACS code2 \and more}
\subclass{65D05 \and 65N30}
\end{abstract}

\section{Introduction}
Analyzing the errors of interpolations on $d$-simplices is an important subject in numerical analysis. It is particularly crucial for finite element error analysis. Let us briefly outline the problems considered in this paper using the Lagrange interpolation operator.

Let $d \in \{ 1,2,3\}$. Let $\widehat{T} \subset \mathbb{R}^d$ and  $T_0 \subset \mathbb{R}^d$ be a reference element and a simplex, respectively, that are affine equivalent. Let us consider two Lagrange finite elements $\{ \widehat{T},\widehat{P} := \mathcal{P}^k, \widehat{\Sigma} \}$ and $\{ {T}_0,{P} := \mathcal{P}^k, {\Sigma} \}$ with associated normed vector spaces $V(\widehat{T}) := \mathcal{C}(\widehat{T})$ and $V(T_0) := \mathcal{C}(T_0)$ with $k \in \mathbb{N}$, where $\mathcal{P}^m$ is the space of polynomials with degree at most $m \in \mathbb{N}_0 := \mathbb{N} \cup \{ 0 \}$. For $\hat{\varphi} \in V(\widehat{T})$, we use the correspondences
\begin{align*}
\displaystyle
&(\varphi_0:T_0 \to \mathbb{R})  \to (\hat{\varphi} := {\varphi}_0 \circ {\Phi}:\widehat{T} \to \mathbb{R}),
\end{align*}
where $\Phi$ is an affine mapping. Let $I_{\widehat{T}}^k: V(\widehat{T}) \to \mathcal{P}^k$ and $I_{{T}_0}^k: V({T}_0) \to \mathcal{P}^k$ be the corresponding Lagrange interpolation operators. Details can be found in Section \ref{sec23=gene}.

We first consider the case in which $d=1$. Let $\Omega := (0,1) \subset \mathbb{R}$. For $N \in \mathbb{N}$, let $\mathbb{T}_h = \{  0 = x_0 \< x_1 \< \cdots \< x_N \< x_{N+1} = 1 \}$ be a mesh of $\overline{\Omega}$ such as
\begin{align*}
\displaystyle
\overline{\Omega} := \bigcup_{i=1}^N T_0^i, \quad \Int T_0^i \cap \Int T_0^j = \emptyset \quad \text{for $i \neq j$},
\end{align*}
where $T_0^i := [x_i , x_{i+1}]$ for $0 \leq i \leq N$. We denote $h_i := x_{i+1} - x_i $ for $0 \leq i \leq N$. If we set $x_j := \frac{j}{N+1}$ for $j=0,1,\ldots,N,N+1$, the mesh  $\mathbb{T}_h$ is said to be the uniform mesh. If we set $x_j := g \left( \frac{j}{N+1} \right)$ for $j=1,\ldots,N,N+1$ with a grading function $g$, the mesh $\mathbb{T}_h$ is said to be the graded mesh with respect to $x=0$; see \cite{BabSur94}. In particular, when $g(y) := y^{\varepsilon}$ ($\varepsilon \> 0$), the mesh is called the radical mesh. To obtain the Lagrange interpolation error estimates, we impose standard assumptions and specify that $\ell, m \in \mathbb{N}_0$ and $p , q \in [1,\infty]$ such that
\begin{align}
\displaystyle
0 \leq m \leq \ell + 1: \quad W^{\ell +1,p}(\widehat{T}) \hookrightarrow W^{m,q} (\widehat{T}). \label{cond=01}
\end{align}
Under these assumptions, the following holds for any $\varphi_0 \in W^{\ell+1,p}(T_0^i)$ with $\hat{\varphi} = \varphi_0 \circ \Phi$:
\begin{align}
\displaystyle
 |\varphi_0 - {I}_{T_0^i}^k \varphi_0 |_{W^{m,q}(T_0^i)}
&\leq c h_{i}^{\frac{1}{q} - \frac{1}{p} + \ell+1-m} | \varphi_0 |_{W^{\ell+1,p}(T_0^i)}. \label{inter=1d}
\end{align}
The proof of this statement is standard; see \cite{ErnGue04}. When $p = q$, it is possible to obtain optimal error estimates even if the scale is different for each element. When $q \> p$, the order of convergence of the interpolation operator may deteriorate.

We now consider the cases in which $d=2,3$. Let $\Omega \subset \mathbb{R}^d$ be a bounded polyhedral domain. Let $\mathbb{T}_h = \{ T_0 \}$ be a simplicial mesh of $\overline{\Omega}$ made up of closed $d$-simplices, such as
\begin{align*}
\displaystyle
\overline{\Omega} = \bigcup_{T_0 \in \mathbb{T}_h} T_0
\end{align*}
with $h := \max_{T_0 \in \mathbb{T}_h} h_{T_0}$, where $ h_{T_0} := \diam(T_0)$. For simplicity, we assume that $\mathbb{T}_h$ is conformal. That is, $\mathbb{T}_h$ is a simplicial mesh of $\overline{\Omega}$ without hanging nodes. Let $\widehat{T} \subset \mathbb{R}^d$ be the reference element defined in Section \ref{sec_2} and $\Phi$ be the affine mapping defined in Eq.~\eqref{affine11}. For any $T_0 \in \mathbb{T}_h$, it holds that $T_0 = \Phi (\widehat{T})$. Under the standard assumptions and Eq.~\eqref{cond=01}, the following holds for any $\varphi_0 \in W^{\ell+1,p}(T_0)$ with $\hat{\varphi} = \varphi_0 \circ \Phi$:
\begin{align}
\displaystyle
 |\varphi_0 - {I}_{T_0}^k \varphi_0 |_{W^{m,q}(T_0)}
&\leq c |T_0|^{\frac{1}{q} - \frac{1}{p}} \left( \frac{\alpha_{\max}}{\alpha_{\min}} \right)^m \left( \frac{H_{T_0}}{h_{T_0}} \right)^m h_{T_0}^{\ell+1-m} | \varphi_0 |_{W^{\ell+1,p}(T_0)}, \label{inter=23d}
\end{align}
where $|T_0|$ is the measure of $T_0$, the parameters $\alpha_{\max}$ and $\alpha_{\min}$ are defined in Eq.~\eqref{cla42}, and the parameter $H_{T_0}$ is as proposed in a recent paper \cite{IshKobTsu}; see Section \ref{sec24} for a definition. The proof of estimate \eqref{inter=23d} can be found in Section \ref{sec_cla}. Compared with the one-dimensional case, the quantities $\alpha_{\max}/\alpha_{\min}$ and $H_{T_0}/h_{T_0}$ negatively affect the order of convergence and do not appear in Eq.\eqref{inter=1d}. The two quantities $\alpha_{\max}/\alpha_{\min}$ and $H_{T_0}/h_{T_0}$ are considered in Section \ref{goodbad}. As a mesh condition, the \textit{shape-regularity condition} is widely used and well known. This condition states that there exists a constant $\gamma \> 0$ such that
\begin{align}
\displaystyle
\rho_{T_0} \geq \gamma h_{T_0} \quad \forall \mathbb{T}_h \in \{ \mathbb{T}_h \}, \quad \forall T_0 \in \mathbb{T}_h, \label{intro2}
\end{align}
where $\rho_{T_0}$ is the radius of the inscribed ball of $T_0$. Under this condition, it holds that
\begin{align}
\displaystyle
 |\varphi_0 - {I}_{T_0}^k \varphi_0 |_{W^{m,q}(T_0)}
&\leq c |T_0|^{\frac{1}{q} - \frac{1}{p}} h_{T_0}^{\ell+1-m} | \varphi_0 |_{W^{\ell+1,p}(T_0)}; \label{intro3}
\end{align}
see Section \ref{Isomesh}. If condition \eqref{intro2} is violated (i.e., the simplex becomes too flat as $h_{T_0} \to 0$), the quantity 
\begin{align*}
\displaystyle
\left( \frac{\alpha_{\max}}{\alpha_{\min}} \right)^m \left( \frac{H_{T_0}}{h_{T_0}} \right)^m h_{T_0}^{\ell+1-m}
\end{align*}
may diverge  even when $p=q$. The effect of the quantity $|T_0|^{\frac{1}{q} - \frac{1}{p}}$ on the interpolation error estimates is considered in Section \ref{qpinflu}.

In some cases, it is not necessary for condition \eqref{intro2} to hold to obtain Eq.~\eqref{intro3}.  The shape-regularity condition can be relaxed to the \textit{maximum-angle condition}, as stated in Eqs.~\eqref{M.A.C.2} and \eqref{M.A.C.3}, for both two-dimensional \cite{BabAzi76} and three-dimensional cases \cite{Kri92}. Anisotropic interpolation theory has also been developed \cite{ApeDob92,Ape99,CheShiZha04}. The idea of Apel \textit{et al.} is to construct a set of functionals satisfying conditions  \eqref{main39}, \eqref{main40}, and \eqref{main41}. The introduction of these functionals makes it possible to remove the quantity $\alpha_{\max}/\alpha_{\min}$. Under the conditions of the maximum angle and coordinate system, anisotropic interpolation error estimates can then be deduced (e.g., see \cite{Ape99}).

In contrast, this paper proposes anisotropic interpolation error estimates using the new parameter under conditions  \eqref{main39}, \eqref{main40}, and \eqref{main41} and Assumption \ref{ass1}; i.e., we derive the following anisotropic error  estimate (Theorem B, in particular, Corollary \ref{coro1}):
\begin{align}
\displaystyle
&| {\varphi}_0 - I_{{T}_0}^k {\varphi}_0 |_{W^{m,q}({T}_0)} \nonumber \\
&\quad \leq  c |T_0|^{\frac{1}{q} - \frac{1}{p}} \left( \frac{H_{T_0}}{h_{T_0}} \right)^m \sum_{|\gamma| = \ell-m} \mathscr{H}^{\gamma}  | \partial^{\gamma}  {(\varphi_0 \circ \Phi_{T_0})} |_{W^{ m ,p}(\Phi_{T_0}^{-1}(T_0))}, \label{intro=06}
\end{align}
where $\Phi_{T_0}$ is defined in Eq.~\eqref{affine10}, $\gamma := (\gamma_1 , \ldots, \gamma_d) \in \mathbb{N}_0^d$ is a multi-index, and $\mathscr{H}^{\gamma}$ is specified in Definition \ref{def=mathscr}. Theorem B applies to interpolations other than the Lagrange interpolation, and the basis for the proof of Theorem B is the \textit{scaling argument} described in Section \ref{sec3=SA}.

Because the new geometric parameter is used in the interpolation error analysis, the coefficient $c$ used in the error estimation is independent of the geometry of the simplices, and the error estimations obtained may therefore be applied to arbitrary meshes, including very ``flat'' or anisotropic simplices. Furthermore, we are naturally able to consider the following geometric condition as being sufficient to obtain optimal order estimates (when $p=q$): there exists $\gamma_0 \> 0$ such that
\begin{align}
\displaystyle
\frac{H_{T_0}}{h_{T_0}} \leq \gamma_0 \quad \forall \mathbb{T}_h \in \{ \mathbb{T}_h \}, \quad \forall T_0 \in \mathbb{T}_h. \label{NewGeo}
\end{align}
Condition \eqref{NewGeo} appears to be simpler than the maximum-angle condition. Furthermore, the quantity $H_{T_0} / h_{T_0}$ can be easily calculated in the numerical process of finite element methods. Therefore, the new condition may be useful. A recent paper \cite{IshKobSuzTsu21} showed that the new condition is satisfied if and only if the maximum-angle condition holds. We expect the new mesh condition to become an alternative to the maximum-angle condition. 

Furthermore, under Assumption \ref{ass1}, component-wise inverse inequalities can be deduced as (see Section \ref{sec_inv}):
\begin{align*}
\displaystyle
\| \partial^{\gamma} \varphi_h \|_{L^q(T)} \leq C^{IVC}  |T|^{\frac{1}{q} - \frac{1}{p}} \mathscr{H}^{- \gamma} \| {\varphi}_h \|_{L^p({T})}.
\end{align*}

In a previous paper \cite{IshKobTsu}, the present authors developed new interpolation error estimations in a general framework and derived Raviart--Thomas interpolations on $d$-simplices. However, the statement of Theorem 2 in \cite{IshKobTsu} includes a mistake. That is, under standard assumptions, the quantity $\alpha_{\max}/\alpha_{\min}$ cannot be removed.  We need to modify the statement of this theorem to correct this error.  The current paper presents Theorems A (see Section \ref{sec_cla}) and B (see Section \ref{sec_main}), which replace Theorem 2 of \cite{IshKobTsu}. In Section \ref{sec_mis}, we explain the inaccuracies in the proof of Theorem 2 in \cite{IshKobTsu} and describe how the results can be recovered using our Theorems A and B. Furthermore, the Babu\v{s}ka and Aziz technique is generally not applicable on anisotropic meshes in the proof of Theorem 3 in \cite{IshKobTsu}. Details will be discussed in a coming paper \cite{Ish21}. 

When there is no ambiguity, we use the notation and definitions given in \cite{IshKobTsu}. Throughout this paper, $c$ denotes a constant independent of $h$ (defined later), unless specified otherwise. These values may change in each context. $\mathbb{R}_+$ is the set of positive real numbers.

\section{Strategy for constructing anisotropic interpolation theory} \label{sec_2}
In standard interpolation theory, one introduces an affine mapping that connects the reference element to the mesh element. However, on anisotropic meshes, the interpolation errors may be overestimated. Therefore, our strategy is to divide the transformation into three affine mappings.

\begin{comment}
\subsection{Mesh}
Let $\Omega \subset \mathbb{R}^d$, $d \in \{ 2 , 3 \}$, be a bounded polyhedral domain. Let $\mathbb{T}_h = \{ T_0 \}$ be a simplicial mesh of $\overline{\Omega}$, made up of closed $d$-simplices, such as
\begin{align*}
\displaystyle
\overline{\Omega} = \bigcup_{T_0 \in \mathbb{T}_h} T_0,
\end{align*}
with $h := \max_{T_0 \in \mathbb{T}_h} h_{T_0}$, where $ h_{T_0} := \diam(T_0)$. For simplicity, we assume that $\mathbb{T}_h$ is conformal. That is, $\mathbb{T}_h$ is a simplicial mesh of $\overline{\Omega}$ without hanging nodes.
\end{comment}

\subsection{Standard positions of simplices} \label{standard}
We recall \cite[Section 3]{IshKobTsu}. Let us first define a diagonal matrix $ \widehat{A}^{(d)}$ as
\begin{align}
\displaystyle
\widehat{A}^{(d)} :=  \diag (\alpha_1,\ldots,\alpha_d), \quad \alpha_i \in \mathbb{R}_+ \quad \forall i.  \label{mesh1}
\end{align}

\subsubsection{Two-dimensional case} \label{reference2d}
Let $\widehat{T} \subset \mathbb{R}^2$ be the reference triangle with vertices $\hat{x}_1 := (0,0)^T$, $\hat{x}_2 := (1,0)^T$, and $\hat{x}_3 := (0,1)^T$. 

Let $\widetilde{\mathfrak{T}}^{(2)}$ be the family of triangles
\begin{align*}
\displaystyle
\widetilde{T} = \widehat{A}^{(2)} (\widehat{T})
\end{align*}
with vertices $\tilde{x}_1 := (0,0)^T$, $\tilde{x}_2 := (\alpha_1,0)^T$, and $\tilde{x}_3 := (0,\alpha_2)^T$.

We next define the regular matrices $\widetilde{A} \in \mathbb{R}^{2 \times 2}$ by
\begin{align}
\displaystyle
\widetilde{A} :=
\begin{pmatrix}
1 & s \\
0 & t \\
\end{pmatrix}, \label{mesh2}
\end{align}
with parameters
\begin{align*}
\displaystyle
s^2 + t^2 = 1, \quad t \> 0.
\end{align*}
For $\widetilde{T} \in \widetilde{\mathfrak{T}}^{(2)}$, let $\mathfrak{T}^{(2)}$ be the family of triangles
\begin{align*}
\displaystyle
T &= \widetilde{A} (\widetilde{T})
\end{align*}
with vertices $x_1 := (0,0)^T, \ x_2 := (\alpha_1,0)^T, \ x_3 :=(\alpha_2 s , \alpha_2 t)^T$. We then have that $\alpha_1 = |x_1 - x_2| \> 0$ and $\alpha_2 = |x_1 - x_3| \> 0$. 

\subsubsection{Three-dimensional case} \label{reference3d}
Let $\widehat{T}_1$ and $\widehat{T}_2$ be reference tetrahedra with the following vertices:
\begin{description}
   \item[(\roman{sone})] $\widehat{T}_1$ has the vertices $\hat{x}_1 := (0,0,0)^T$, $\hat{x}_2 := (1,0,0)^T$, $\hat{x}_3 := (0,1,0)^T$, $\hat{x}_4 := (0,0,1)^T$;
 \item[(\roman{stwo})] $\widehat{T}_2$ has the vertices $\hat{x}_1 := (0,0,0)^T$, $\hat{x}_2 := (1,0,0)^T$, $\hat{x}_3 := (1,1,0)^T$, $\hat{x}_4 := (0,0,1)^T$.
\end{description}

Let $\widetilde{\mathfrak{T}}_i^{(3)}$, $i=1,2$, be the family of triangles
\begin{align*}
\displaystyle
\widetilde{T}_i = \widehat{A}^{(3)} (\widehat{T}_i), \quad i=1,2,
\end{align*}
with vertices 
\begin{description}
   \item[(\roman{sone})] $\tilde{x}_1 := (0,0,0)^T$, $\tilde{x}_2 := (\alpha_1,0,0)^T$, $\tilde{x}_3 := (0,\alpha_2,0)^T$, and $\tilde{x}_4 := (0,0,\alpha_3)^T$;
 \item[(\roman{stwo})]  $\tilde{x}_1 := (0,0,0)^T$, $\tilde{x}_2 := (\alpha_1,0,0)^T$, $\tilde{x}_3 := (\alpha_1,\alpha_2,0)^T$, and $\tilde{x}_4 := (0,0,\alpha_3)^T$.
\end{description}

We next define the regular matrices $\widetilde{A}_1, \widetilde{A}_2 \in \mathbb{R}^{3 \times 3}$ by
\begin{align}
\displaystyle
\widetilde{A}_1 :=
\begin{pmatrix}
1 & s_1 & s_{21} \\
0 & t_1  & s_{22}\\
0 & 0  & t_2\\
\end{pmatrix}, \
\widetilde{A}_2 :=
\begin{pmatrix}
1 & - s_1 & s_{21} \\
0 & t_1  & s_{22}\\
0 & 0  & t_2\\
\end{pmatrix} \label{mesh3}
\end{align}
with parameters
\begin{align*}
\displaystyle
\begin{cases}
s_1^2 + t_1^2 = 1, \ s_1 \> 0, \ t_1 \> 0, \ \alpha_2 s_1 \leq \alpha_1 / 2, \\
s_{21}^2 + s_{22}^2 + t_2^2 = 1, \ t_2 \> 0, \ \alpha_3 s_{21} \leq \alpha_1 / 2.
\end{cases}
\end{align*}
For $\widetilde{T}_i \in \widetilde{\mathfrak{T}}_i^{(3)}$, $i=1,2$, let $\mathfrak{T}_i^{(3)}$, $i=1,2$, be the family of tetrahedra
\begin{align*}
\displaystyle
T_i &= \widetilde{A}_i (\widetilde{T}_i), \quad i=1,2,
\end{align*}
with vertices
\begin{align*}
\displaystyle
&x_1 := (0,0,0)^T, \ x_2 := (\alpha_1,0,0)^T, \ x_4 := (\alpha_3 s_{21}, \alpha_3 s_{22}, \alpha_3 t_2)^T, \\
&\begin{cases}
x_3 := (\alpha_2 s_1 , \alpha_2 t_1 , 0)^T \quad \text{for case (\roman{sone})}, \\
x_3 := (\alpha_1 - \alpha_2 s_1, \alpha_2 t_1,0)^T \quad \text{for case (\roman{stwo})}.
\end{cases}
\end{align*}
We then have $\alpha_1 = |x_1 - x_2| \> 0$, $\alpha_3 = |x_1 - x_4| \> 0$, and
\begin{align*}
\displaystyle
\alpha_2 =
\begin{cases}
|x_1 - x_3| \> 0  \quad \text{for case (\roman{sone})}, \\
|x_2 - x_3| \> 0  \quad \text{for case (\roman{stwo})}.
\end{cases}
\end{align*}

In the following, we impose conditions for $T \in \mathfrak{T}^{(2)}$ in the two-dimensional case and $T \in \mathfrak{T}_1^{(3)} \cup \mathfrak{T}_2^{(3)} =: \mathfrak{T}^{(3)}$ in the three-dimensional case. 

\begin{cond}[Case in which $d=2$] \label{cond1}
Let $T \in \mathfrak{T}^{(2)}$ with the vertices $x_i$ ($i=1,\ldots,3$) introduced in Section \ref{reference2d}. We assume that $\overline{x_2 x_3}$ is the longest edge of $T$; i.e., $ h_T := |x_2 - x_ 3|$. Recall that $\alpha_1 = |x_1 - x_2|$ and $\alpha_2 = |x_1 - x_3|$. We then assume that $\alpha_2 \leq \alpha_1$. Note that $\alpha_1 = \mathcal{O}(h_T)$. 
\end{cond}
\begin{cond}[Case in which $d=3$] \label{cond2}
Let $T \in \mathfrak{T}^{(3)}$ with the vertices $x_i$ ($i=1,\ldots,4$) introduced in Section \ref{reference3d}. Let $L_i$ ($1 \leq i \leq 6$) be the edges of $T$. We denote by $L_{\min}$  the edge of $T$ that has the minimum length; i.e., $|L_{\min}| = \min_{1 \leq i \leq 6} |L_i|$. We set $\alpha_2 := |L_{\min}|$ and assume that 
\begin{align*}
\displaystyle
&\text{the end points of $L_{\min}$ are either $\{ x_1 , x_3\}$ or $\{ x_2 , x_3\}$}.
\end{align*}
Among the four edges that share an end point with $L_{\min}$, we take the longest edge $L^{({\min})}_{\max}$. Let $x_1$ and $x_2$ be the end points of edge $L^{({\min})}_{\max}$. Thus, we have that
\begin{align*}
\displaystyle
\alpha_1 = |L^{(\min)}_{\max}| = |x_1 - x_2|.
\end{align*}
Consider cutting $\mathbb{R}^3$ with the plane that contains the midpoint of edge $L^{(\min)}_{\max}$ and is perpendicular to the vector $x_1 - x_2$. We then have two cases: 
\begin{description}
  \item[(Type \roman{sone})] $x_3$ and $x_4$  belong to the same half-space;
  \item[(Type \roman{stwo})] $x_3$ and $x_4$  belong to different half-spaces.
\end{description}
In each respective case, we set
\begin{description}
  \item[(Type \roman{sone})] $x_1$ and $x_3$ as the end points of $L_{\min}$, that is, $\alpha_2 =  |x_1 - x_3| $;
  \item[(Type \roman{stwo})] $x_2$ and $x_3$ as the end points of $L_{\min}$, that is, $\alpha_2 =  |x_2 - x_3| $.
\end{description}
Finally, we set $\alpha_3 = |x_1 - x_4|$. Note that we implicitly assume that $x_1$ and $x_4$ belong to the same half-space. In addition, note that $\alpha_1 = \mathcal{O}(h_T)$.
\end{cond}

\begin{comment}
Let $\Omega \subset \mathbb{R}^d$, $d \in \{ 2 , 3 \}$, be a bounded polyhedral domain. Let $\mathbb{T}_h = \{ T_0 \}$ be a simplicial mesh of $\overline{\Omega}$, made up of closed $d$-simplices, such as
\begin{align*}
\displaystyle
\overline{\Omega} = \bigcup_{T_0 \in \mathbb{T}_h} T_0,
\end{align*}
with $h := \max_{T_0 \in \mathbb{T}_h} h_{T_0}$, where $ h_{T_0} := \diam(T_0)$. For simplicity, we assume that $\mathbb{T}_h$ is conformal. That is, $\mathbb{T}_h$ is a simplicial mesh of $\overline{\Omega}$ without hanging nodes.
\end{comment}

\begin{Rem} \label{Rem11111}
%We impose Condition \ref{Cond334}. 
Let $\mathbb{T}_h$ be a conformal mesh. We assume that any simplex $T_0 \in \mathbb{T}_h$ is transformed into $T_1 \in \mathfrak{T}^{(2)}$ such that Condition \ref{cond1} is satisfied (in the two-dimensional case) or $T_i \in \mathfrak{T}_i^{(3)}$, $i=1,2$, such that Condition \ref{cond2} is satisfied (in the three-dimensional case) through appropriate rotation, translation, and mirror imaging. Note that none of the lengths of the edges of a simplex or the measure of the simplex is changed by the transformation.
\end{Rem}

\begin{ass} \label{ass1}
In anisotropic interpolation error analysis, we may impose the following geometric conditions for the simplex $T$:
\begin{enumerate}
 \item If $d=2$, there are no additional conditions;
 \item If $d=3$, there exists a positive constant $M$, independent of $h_T$, such that $|s_{22}| \leq M \frac{\alpha_2 t_1}{\alpha_3}$. Note that if $s_{22} \neq 0$, this condition means that the order with respect to $h_T$ of $\alpha_3$ coincides with the order of $\alpha_2$, whereas if $s_{22} = 0$, the order of $\alpha_3$ may be different from that of $\alpha_2$. 
\end{enumerate}
\end{ass}

\subsection{Affine mappings}  \label{APsubsection34}
In our strategy, we adopt the following affine mappings.

\begin{defi}[Affine mappings] \label{affine=defi}
Let $\widetilde{T}, \widehat{T} \subset \mathbb{R}^d$ be the simplices defined in Sections \ref{reference2d} and \ref{reference3d}. That is, 
\begin{align*}
\displaystyle
\widetilde{T} = \widehat{\Phi} (\widehat{T}), \quad {T} = \widetilde{\Phi} (\widetilde{T}) \quad  \text{with} \quad \tilde{x} := \widehat{\Phi}(\hat{x}) := \widehat{{A}}^{(d)} \hat{x}, \quad x := \widetilde{\Phi}(\tilde{x}) := \widetilde {{A}} \tilde{x}.
\end{align*}
We then define an affine mapping $\Phi_T: \widehat{T} \to T$ as
\begin{align}
\displaystyle
\Phi_T := \widetilde{\Phi} \circ \widehat{\Phi}: \widehat{T} \to T, \  {x} :={\Phi}_T(\hat{x}) := {{A}}_T \hat{x}, \quad {A}_T := \widetilde {{A}} \widehat{{A}}^{(d)}. \label{affine9}
\end{align}
Furthermore, let $\Phi_{T_0}$ be an affine mapping defined as
\begin{align}
\displaystyle
\Phi_{T_0}: T \ni x \mapsto {A}_{T_0} x + b_{T_0} \in T_0, \label{affine10}
\end{align}
where $b_{T_0} \in \mathbb{R}^d$ and ${A}_{T_0} \in O(d)$ is a rotation and mirror imaging matrix. We then define an affine mapping $\Phi: \widehat{T} \to T_0$ as
\begin{align}
\displaystyle
\Phi := {\Phi}_{T_0} \circ {\Phi}_T: \widehat{T} \to T_0, \ x^{(0)} := \Phi (\hat{x}) =  ({\Phi}_{T_0} \circ {\Phi}_T)(\hat{x}) = {A} \hat{x} + b_{T_0}, \label{affine11}
\end{align}
where ${A} := {A}_{T_0} {A}_T$.
\end{defi}

\subsection{Finite element generation} \label{sec23=gene}
We follow the procedure described in \cite[Section 1.4.1 and 1.2.1]{ErnGue04}; see also \cite[Section 3.5]{IshKobTsu}.

For the reference element $\widehat{T}$ defined in Sections \ref{reference2d} and \ref{reference3d}, let $\{ \widehat{T} , \widehat{{P}} , \widehat{\Sigma} \}$ be a fixed reference finite element, where $\widehat{{P}} $ is a vector space of functions $\hat{p}: \widehat{T} \to \mathbb{R}^n$ for some positive integer $n$ (typically $n=1$ or $n=d$) and $\widehat{\Sigma}$ is a set of $n_{0}$ linear forms $\{ \hat{\chi}_1 , \ldots , \hat{\chi}_{n_0} \}$ such that
\begin{align*}
\displaystyle
\widehat{{P}} \ni \hat{p} \mapsto (\hat{\chi}_1(\hat{p}) , \ldots , \hat{\chi}_{n_0}(\hat{p}))^T \in \mathbb{R}^{n_0}
\end{align*}
is bijective; i.e., $\widehat{\Sigma}$ is a basis for $\mathcal{L}(\widehat{P};\mathbb{R})$. Further, we denote by $\{ \hat{\theta}_1 , \ldots, \hat{\theta}_{n_0} \}$ in $\widehat{{P}}$ the local ($\mathbb{R}^n$-valued) shape functions such that
\begin{align*}
\displaystyle
\hat{\chi}_i(\hat{\theta}_j) = \delta_{ij}, \quad 1 \leq i,j \leq n_0.
\end{align*}
Let $V(\widehat{T})$ be a normed vector space of functions $\hat{\varphi}: \widehat{T} \to \mathbb{R}^n$ such that $\widehat{P} \subset V(\widehat{T})$ and the linear forms $\{ \hat{\chi}_1 , \ldots , \hat{\chi}_{n_0} \}$ can be extended to $V(\widehat{T})^{\prime}$. The local interpolation operator ${I}_{\widehat{T}}$ is then defined by
\begin{align}
\displaystyle
{I}_{\widehat{T}} : V(\widehat{T}) \ni \hat{\varphi} \mapsto \sum_{i=1}^{n_0} \hat{\chi}_i (\hat{\varphi}) \hat{\theta}_i \in \widehat{{P}}. \label{int1}
\end{align}
It is obvious that
\begin{align}
\displaystyle
& \hat{\chi}_i ({I}_{\widehat{T}} \hat{\varphi}) =  \hat{\chi}_i (\hat{\varphi}) \quad \forall \hat{\varphi} \in V(\widehat{T}), \quad i=1,\ldots,n_0, \label{int512} \\
& I_{\widehat{T}} \hat{p} = \hat{p} \quad \forall \hat{p} \in \widehat{P}. \label{int513}
\end{align}
Let ${\Phi}$ be the affine mapping defined in Eq.~\eqref{affine11}. For $T_0  = {\Phi} (\widehat{T})$, we first define a Banach space $V(T_0)$ of $\mathbb{R}^n$-valued functions that is the counterpart of $V(\widehat{T})$ and define a linear bijection mapping by
\begin{align*}
\displaystyle
\psi  := \psi_{\widehat{T}} \circ \psi_{\widetilde{T}} \circ \psi_{T} : V(T_0) \ni \varphi \mapsto \hat{\varphi} := \psi(\varphi) := \varphi \circ \Phi \in   V(\widehat{T})
\end{align*}
with the three linear bijection mappings
\begin{align*}
\displaystyle
&\psi_{{T}} : V({T}_0) \ni \varphi_0 \mapsto {\varphi} := \psi_{{T}}(\varphi_0) := \varphi_0 \circ {\Phi}_{T_0} \in   V({T}), \\
&\psi_{\widetilde{T}} : V({T}) \ni \varphi \mapsto \tilde{\varphi} := \psi_{\widetilde{T}}(\varphi) := \varphi \circ \widetilde{\Phi} \in   V(\widetilde{T}), \\
&\psi_{\widehat{T}} : V(\widetilde{T}) \ni \tilde{\varphi} \mapsto \hat{\varphi} :=  \psi_{\widehat{T}} (\tilde{\varphi}) :=  \tilde{\varphi} \circ \widehat{\Phi} \in V(\widehat{T}).
\end{align*}
The triple $\{ \widetilde{T} , \widetilde{P} , \widetilde{\Sigma} \}$ is defined as
\begin{align*}
\displaystyle
\begin{cases}
\displaystyle
\widetilde{T} = \widehat{\Phi}(\widehat{T}); \\
\displaystyle
 \widetilde{P} = \{ \psi_{\widehat{T}}^{-1}(\hat{p}) ; \ \hat{p} \in \widehat{{P}}\}; \\
\displaystyle
\widetilde{\Sigma} = \{ \{ \tilde{\chi}_{i} \}_{1 \leq i \leq n_0}; \ \tilde{\chi}_{i} = \hat{\chi}_i(\psi_{\widehat{T}}(\tilde{p})), \forall \tilde{p} \in \widetilde{P}, \hat{\chi}_i \in \widehat{\Sigma} \}.
\end{cases}
\end{align*}
The triples  $\{ {T} , {P} , {\Sigma} \}$ and  $\{ {T}_0 , {P}_0 ,  {\Sigma}_0 \}$ are similarly defined. These triples are finite elements and the local shape functions are $\tilde{\theta}_{i} = \psi_{\widehat{T}}^{-1}(\hat{\theta}_i)$, $\theta_{i} = \psi_{\widetilde{T}}^{-1}(\tilde{\theta}_i)$, and $\theta_{0,i} := \psi_{T}^{-1} (\theta_i)$ for $1 \leq i \leq n_0$, and the associated local interpolation operators are respectively defined by
\begin{align}
\displaystyle
{I}_{\widetilde{T}} : V(\widetilde{T}) \ni \tilde{\varphi} \mapsto {I}_{\widetilde{T}} \tilde{\varphi} &:= \sum_{i=1}^{n_0} \tilde{\chi}_{i}(\tilde{\varphi}) \tilde{\theta}_{i} \in \widetilde{P}, \label{int2} \\
{I}_{T} : V(T) \ni \varphi \mapsto {I}_{T} \varphi &:= \sum_{i=1}^{n_0} \chi_{i}(\varphi) \theta_{i} \in {P}, \label{int3} \\
{I}_{T_0} : V(T_0) \ni \varphi_0 \mapsto {I}_{T_0} \varphi_0 &:= \sum_{i=1}^{n_0} \chi_{0,i}(\varphi_0) \theta_{i} \in {P}_0. \label{int4}
\end{align}

\begin{prop} \label{prop512}
The diagrams
\begin{align*}
\displaystyle
\xymatrix{
V(T_0)\ar[r]^-{\psi_{{T}}}\ar[d]_-{I_{T_0}}\ar@{}|{}&V(T)\ar[r]^-{\psi_{\widetilde{T}}}\ar[d]_-{I_{T}}\ar@{}|{}&V(\widetilde{T})\ar[r]^-{\psi_{\widehat{T}}}\ar[d]_-{{I}_{\widetilde{T}}}\ar@{}|{}&V(\widehat{T})\ar[d]^-{{I}_{\widehat{T}}} \\
P_0\ar[r]_-{\psi_{{T}}}&P\ar[r]_-{\psi_{\widetilde{T}}}&\widetilde{{P}}\ar[r]_-{\psi_{\widehat{T}}}&\widehat{{P}}
}
\end{align*}
commute.
\end{prop}

\begin{pf*}
See, for example, \cite[Proposition 1.62]{ErnGue04}.
\qed
\end{pf*}

\subsection{New parameters} \label{sec24}
 In a previous paper \cite{IshKobTsu}, we proposed two geometric parameters, 
 \begin{defi} \label{defi1}
 The parameter $H_T$ is defined as
\begin{align*}
\displaystyle
H_T := \frac{\prod_{i=1}^d \alpha_i}{|T|} h_T,
\end{align*}
and the parameter $H_{T_0}$ is defined as
\begin{align*}
\displaystyle
H_{T_0} := \frac{h_{T_0}^2}{|T_0|} \min_{1 \leq i \leq 3} |L_i|  \quad \text{if $d=2$}, \quad H_{T_0} := \frac{h_{T_0}^2}{|T_0|} \min_{1 \leq i , j \leq 6, i \neq j} |L_i| |L_j| \quad \text{if $d=3$}
\end{align*}
where $L_i$ denotes the edges of the simplex $T_0 \subset \mathbb{R}^d$. 
\end{defi}

The following lemma shows the equivalence between $H_T$ and $H_{T_0}$.

 \begin{lem} \label{lem2}
 It holds that
\begin{align*}
\displaystyle
\frac{1}{2} H_{T_0} \< H_T  \< 2 H_{T_0}.
\end{align*}
Furthermore, in the two-dimensional case, $H_{T_0}$ is equivalent to the circumradius $R_2$ of $T_0$.
 \end{lem}
 
 \begin{pf*}
The proof can be found in \cite[Lemma 3]{IshKobTsu}.
\qed
\end{pf*}

\begin{Rem}
We set
\begin{align*}
\displaystyle
H(h) := \max_{T_0 \in \mathbb{T}_h} H_{T_0}.
\end{align*}
As we stated in the Introduction, if the maximum-angle condition is violated, the parameter $H(h)$ may diverge as $h \to 0$ on anisotropic meshes. Therefore, imposing the maximum-angle condition for mesh partitions guarantees the convergence of finite element methods \cite{Apeet21}. Reference \cite{BabAzi76} studied cases in which the finite element solution may not converge to the exact solution.	
\end{Rem}

We now state the following theorem concerning the new condition.

\begin{thr} \label{thr3}
Condition \eqref{NewGeo} holds if and only if there exist $0 \< \gamma_1, \gamma_2 \< \pi$ such that 
\begin{align}
\displaystyle
d=2: \quad \theta_{T_0,\max} \leq \gamma_1 \quad \forall \mathbb{T}_h \in \{ \mathbb{T}_h \}, \quad \forall T_0 \in \mathbb{T}_h, \label{M.A.C.2}
\end{align}
where $\theta_{T_0,\max}$ is the maximum angle of $T_0$, and 
\begin{align}
\displaystyle
d=3: \quad \theta_{T_0, \max} \leq \gamma_2, \quad \psi_{T_0, \max} \leq \gamma_2 \quad \forall \mathbb{T}_h \in \{ \mathbb{T}_h \}, \quad \forall T_0 \in \mathbb{T}_h,  \label{M.A.C.3}
\end{align}
where $\theta_{T_0, \max}$ is the maximum angle of all triangular faces of the tetrahedron $T_0$ and $\psi_{T_0, \max}$ is the maximum dihedral angle of $T_0$. Conditions \eqref{M.A.C.2} and \eqref{M.A.C.3} together constitute the \textit{maximum-angle condition}.
\end{thr}

\begin{pf*}
In the case of $d=2$, we use the previous result presented in \cite{Kri91}; i.e., there exists a constant $\gamma_3  \> 0$ such that
\begin{align*}
\displaystyle
\frac{R_2}{h_{T_0}} \leq \gamma_3 \quad \forall \mathbb{T}_h \in \{ \mathbb{T}_h \}, \quad \forall T_0 \in \mathbb{T}_h,
\end{align*}
if and only if condition \eqref{M.A.C.2} is satisfied. Combining this result with $H_{T_0}$ being equivalent to the circumradius $R_2$ of $T_0$, we have the desired conclusion. In the case of $d=3$, the proof can be found in a recent paper \cite{IshKobSuzTsu21}.
\qed	
\end{pf*}

%\subsection{Euclidean condition number}

\begin{lem} \label{lem351}
It holds that
\begin{subequations} \label{CN331}
\begin{align}
\displaystyle
\| \widehat{{A}}^{(d)} \|_2 &\leq  h_{T}, \quad \| \widehat{{A}}^{(d)} \|_2 \| (\widehat{{A}}^{(d)})^{-1} \|_2 = \frac{\max \{\alpha_1 , \cdots, \alpha_d \}}{\min \{\alpha_1 , \cdots, \alpha_d \}}, \label{CN331a} \\
\| \widetilde{{A}} \|_2 &\leq 
\begin{cases}
\sqrt{2} \quad \text{if $d=2$}, \\
2  \quad \text{if $d=3$},
\end{cases}
\quad \| \widetilde{{A}} \|_2 \| \widetilde{{A}}^{-1} \|_2 \leq
\begin{cases}
\frac{\alpha_1 \alpha_2}{|T|} = \frac{H_{T}}{h_{T}} \quad \text{if $d=2$}, \\
\frac{2}{3} \frac{\alpha_1 \alpha_2 \alpha_3}{|T|} = \frac{2}{3} \frac{H_{T}}{h_{T}} \quad \text{if $d=3$},
\end{cases} \label{CN331b} \\
\| A_{T_0} \|_2 &= 1, \quad \| A_{T_0}^{-1} \|_2 = 1. \label{CN331c}
\end{align}
\end{subequations}
Furthermore, we have
\begin{align}
\displaystyle
| \det ({A}_T) | = | \det(\widetilde{{A}}) | | \det (\widehat{{A}}^{(d)}) | = d ! |T |, \quad | \det ({A}_{T_0}) | = 1. \label{CN332}
\end{align}
\end{lem}

\begin{pf*}
The proof of  \eqref{CN331b} can be found in \cite[(4.4), (4.5), (4.6), and (4.7)]{IshKobTsu}. The inequality \eqref{CN331a} is easily proved.  Because $A_{T_0} \in O(d)$, one easily finds that $A_{T_0}^{-1} \in O(d)$ and recovers Eq.~\eqref{CN331c}. The proof of equality \eqref{CN332} is standard.
\qed
\end{pf*}

For matrix $A \in \mathbb{R}^{d \times d}$, we denote by $[{A}]_{ij}$ the $(i,j)$-component of $A$. We set $\| A \|_{\max} := \max_{1 \leq i,j \leq d} | [{A}]_{ij} |$. Furthermore, we use the inequality
\begin{align}
\displaystyle
 \| {A} \|_{\max} \leq \| {A} \|_2. \label{matrix22}
 \end{align}

\section{Scaling argument} \label{sec3=SA}
This section gives estimates related to a scaling argument corresponding to \cite[Lemma 1.101]{ErnGue04}. The estimates play major roles in our analysis. Furthermore, we use the following inequality (see \cite[Exercise 1.20]{ErnGue04}). Let $0 \< r \leq s$ and $a_i \geq 0$, $i=1,2,\ldots,n$ ($n \in \mathbb{N}$), be real numbers. Then, we have that
\begin{align}
\displaystyle
\left( \sum_{i=1}^n a_i^s \right)^{1/s} \leq \left( \sum_{i=1}^n a_i^r \right)^{1/r}. \label{jensen}
\end{align}	

\begin{lem} \label{int=lem3}
Let $s \geq 0$ and $1 \leq p \leq \infty$. There exist positive constants $c_1$ and $c_2$ such that, for all $T_0 \in \mathbb{T}_h$ and $\varphi_0 \in W^{s,p}(T_0)$,
\begin{align}
\displaystyle
c_1  |\varphi_0 |_{W^{s,p}(T_0)} \leq |\varphi|_{W^{s,p}(T)} &\leq c_ 2 | \varphi_0 |_{W^{s,p}(T_0)} \label{sa22}
\end{align}
with $\varphi = \varphi_0 \circ {\Phi}_{T_0}$. 
\end{lem}

\begin{pf*}
The following inequalities can be found in \cite[Lemma 1.101]{ErnGue04}. There exists a positive constant $c$ such that, for all $T_0 \in \mathbb{T}_h$ and $\varphi_0 \in W^{s,p}(T_0)$,
\begin{align}
\displaystyle
|\varphi|_{W^{s,p}(T)} &\leq c \| {A}_{T_0} \|_2^s |\det ({A}_{T_0})|^{ - \frac{1}{p}} |\varphi_0 |_{W^{s,p}(T_0)}, \label{sa23}\\
|\varphi_0 |_{W^{s,p}(T_0)} &\leq c \| {A}_{T_0}^{-1} \|_2^s |\det ({A}_{T_0})|^{ \frac{1}{p}} |\varphi|_{W^{s,p}(T)} \label{sa24}
\end{align}
with $\varphi = \varphi_0 \circ {\Phi}_{T_0}$. Using Eqs.~\eqref{sa23} and \eqref{sa24} together with Eqs.~\eqref{CN331c} and \eqref{CN332} yields Eq.~\eqref{sa22}.
\qed
\end{pf*}

\begin{lem} \label{int=lem4}
Let $m \in \mathbb{N}_0$ and $p \in [0,\infty)$. Let  $\beta := (\beta_1,\ldots,\beta_d) \in \mathbb{N}_0^d$  be a multi-index with $|\beta| = m$. Then, for any $\hat{\varphi} \in W^{m,p}(\widehat{T})$ with $\tilde{\varphi} = \hat{\varphi} \circ \widehat{\Phi}^{-1}$ and ${\varphi} = \tilde{\varphi} \circ \widetilde{\Phi}^{-1}$, it holds that
\begin{align}
\displaystyle
|{\varphi}|_{W^{m,p}({T})}
&\leq C_1^{SA} |\det(A_T)|^{\frac{1}{p}} \| \widetilde{A}^{-1} \|^m_2 \left( \sum_{|\beta| = m}  ( \alpha^{- \beta} )^p \|  \partial^{\beta} \hat{\varphi} \|^p_{L^p(\widehat{T})} \right)^{1/p}, \label{sa25}
\end{align}
where $C_1^{SA}$ is a constant that is independent of $T$ and $\widetilde{T}$. When $p = \infty$, for any $\hat{\varphi} \in W^{m,\infty}(\widehat{T})$ with $\tilde{\varphi} = \hat{\varphi} \circ \widehat{\Phi}^{-1}$ and ${\varphi} = \tilde{\varphi} \circ \widetilde{\Phi}^{-1}$,  it holds that
\begin{align}
\displaystyle
|{\varphi}|_{W^{m,\infty}({T})} 
%&\leq C_1^{SA,\infty}  \| \widetilde{A}^{-1} \|^m_2 \sum_{|\beta| = m}   \alpha^{- \beta}  \|  \partial^{\beta} \hat{\varphi} \|_{L^{\infty}(\widehat{T})}, \label{sa=infty}
&\leq C_1^{SA,\infty}  \| \widetilde{A}^{-1} \|^m_2 \max_{|\beta| = m} \left(  \alpha^{- \beta}  \|  \partial^{\beta} \hat{\varphi} \|_{L^{\infty}(\widehat{T})} \right), \label{sa=infty}
\end{align}
where $C_1^{SA,\infty}$ is a constant that is independent of $T$ and $\widetilde{T}$. 
\end{lem}

\begin{pf*}
Let $p \in [1,\infty)$. Because the space $\mathcal{C}^{m}(\widehat{T})$ is dense in the space ${W}^{m,p}(\widehat{T})$, we show that Eq.~\eqref{sa25} holds for $\hat{\varphi} \in \mathcal{C}^{m}(\widehat{T})$ with $\tilde{\varphi} = \hat{\varphi} \circ \widehat{\Phi}^{-1}$ and ${\varphi} = \tilde{\varphi} \circ \widetilde{\Phi}^{-1}$. From $\hat{x}_j = \alpha_j^{-1} \tilde{x}_j$, we have that, for any multi-index $\beta$,
\begin{align}
\displaystyle
\partial^{\beta} \tilde{\varphi} &=  \alpha_{1}^{- \beta_1} \cdots \alpha_{d}^{- \beta_d} \partial^{\beta} \hat{\varphi} = \alpha^{- \beta} \partial^{\beta} \hat{\varphi}. \label{hattilde}
\end{align}
Through a change of variable, we obtain
\begin{align}
\displaystyle
|\tilde{\varphi}|_{W^{m,p}(\widetilde{T})}^p
&= \sum_{|\beta| = m}  \| \partial^{\beta} \tilde{\varphi} \|^p_{L^p(\widetilde{T})} 
= | \det (\widehat{A}^{(d)}) | \sum_{|\beta| = m}   (\alpha^{- \beta } )^p \|  \partial^{\beta} \hat{\varphi} \|^p_{L^p(\widehat{T})}. \label{sa27}
\end{align}
From the standard estimate in \cite[Lemma 1.101]{ErnGue04}, we have
\begin{align}
\displaystyle
|{\varphi}|_{W^{m,p}({T})}
&\leq C_1^{SA} |\det(\widetilde{A})|^{\frac{1}{p}} \| \widetilde{A}^{-1} \|^m_2  |\tilde{\varphi}|_{W^{m,p}(\widetilde{T})}. \label{sa28}
\end{align}
Inequality \eqref{sa25} follows from Eqs.~\eqref{sa27} and \eqref{sa28} with Eq.~\eqref{CN332}.

We consider the case in which $p = \infty$. A function $\hat{\varphi} \in W^{m,\infty}(\widehat{T})$ belongs to the space $W^{m,p}(\widehat{T})$ for any $p \in [1,\infty)$.  Therefore, it holds that $\tilde{\varphi} \in W^{m,p}(\widetilde{T})$ for any $p \in [1,\infty)$ and, from Eq.~\eqref{jensen}, we obtain
\begin{align}
\displaystyle
\| \partial^{\gamma} \tilde{ \varphi } \|_{L^{p}(\widetilde{T})} 
&\leq | \tilde{\varphi}|_{W^{|\gamma|,p}(\widetilde{T})} \notag \\
&=  |\det( \widehat{A}^{(d)})|^{\frac{1}{p}} \left( \sum_{|\beta| = |\gamma|}  ( \alpha^{- \beta} )^p \|  \partial^{\beta} \hat{\varphi} \|^p_{L^p(\widehat{T})} \right)^{1/p} \notag \\
&\leq \left( \sup_{1 \leq p} |\det(\widehat{A}^{(d)})|^{\frac{1}{p}} \right) \sum_{|\beta| = |\gamma|}   \alpha^{- \beta}  \|  \partial^{\beta} \hat{\varphi} \|_{L^p(\widehat{T})} \notag\\
&\leq c \left( \sup_{1 \leq p} |\det(\widehat{A}^{(d)})|^{\frac{1}{p}} \right) \sum_{|\beta| = |\gamma|}   \alpha^{- \beta}  \|  \partial^{\beta} \hat{\varphi} \|_{L^{\infty}(\widehat{T})} \< \infty \label{infty31}
\end{align}
for the multi-index $\gamma \in \mathbb{N}_0^d$ with $|\gamma| \leq m$. This implies that the function $\partial^{\gamma} \tilde{\varphi}$ is in the space $L^{\infty}(\widetilde{T})$ for each $|\gamma| \leq m$. Therefore, we have that $\tilde{\varphi} \in W^{m,\infty}(\widetilde{T})$. Taking the limit $p \to \infty$ in Eq.~\eqref{infty31} and using $\lim_{p \to \infty} \| \cdot \|_{L^p(\widetilde{T})} = \| \cdot\|_{L^{\infty}(\widetilde{T})}$, we have
\begin{align}
\displaystyle
| \tilde{\varphi} |_{W^{m,\infty}(\widetilde{T})} \leq c \max_{|\beta| = m} \left(  \alpha^{- \beta}  \|  \partial^{\beta} \hat{\varphi} \|_{L^{\infty}(\widehat{T})} \right). \label{infty35}
\end{align}
From the standard estimate in \cite[Lemma 1.101]{ErnGue04}, we have
\begin{align}
\displaystyle
|{\varphi}|_{W^{m,\infty}({T})}
&\leq c \| \widetilde{A}^{-1} \|^m_2  |\tilde{\varphi}|_{W^{m,\infty}(\widetilde{T})}. \label{sa28inf}
\end{align}
Inequality \eqref{sa=infty} follows from Eqs.~\eqref{infty35} and \eqref{sa28inf}.
\qed
\end{pf*}

We now introduce the following new notation.

\begin{defi} \label{def=mathscr}
We define a parameter $\mathscr{H}_i$, $i=1,\ldots,d$, as
\begin{align*}
\displaystyle
\begin{cases}
\mathscr{H}_1 := \alpha_1, \quad \mathscr{H}_2 := \alpha_2 t \quad \text{if $d=2$}, \\
\mathscr{H}_1 := \alpha_1, \quad \mathscr{H}_2 := \alpha_2 t_1, \quad \mathscr{H}_3 := \alpha_3 t_2 \quad \text{if $d=3$}.
\end{cases}
\end{align*}
For a multi-index $\beta = (\beta_1,\ldots,\beta_d) \in \mathbb{N}_0^d$, we use the notation
\begin{align*}
\displaystyle
\mathscr{H}^{\beta} := \mathscr{H}_1^{\beta_1} \cdots \mathscr{H}_d^{\beta_d}, \quad \mathscr{H}^{- \beta} := \mathscr{H}_1^{- \beta_1} \cdots \mathscr{H}_d^{- \beta_d}.
\end{align*}
We also define  $\alpha^{\beta} :=  \alpha_{1}^{\beta_1} \cdots \alpha_{d}^{\beta_d}$ and $\alpha^{- \beta} :=  \alpha_{1}^{- \beta_1} \cdots \alpha_{d}^{- \beta_d}$. 
\end{defi}

\begin{defi} \label{defi332}
We define vectors $r_n \in \mathbb{R}^d$, $n=1,\ldots,d$, as follows. If $d=2$,
\begin{align*}
\displaystyle
r_1 := (1 , 0)^T, \quad r_2 := (s,t)^T,
\end{align*}
and if $d=3$,
\begin{align*}
\displaystyle
&r_1 := (1 , 0,0)^T, \quad r_3 := ( s_{21}, s_{22} , t_2)^T, \\
&\begin{cases}
r_2 := ( s_1 ,  t_1 , 0)^T \quad \text{for case (\roman{sone})}, \\
r_2 := (- s_1,  t_1,0)^T \quad \text{for case (\roman{stwo})}.
\end{cases}
\end{align*}
Furthermore, we define a directional derivative as
\begin{align*}
\displaystyle
\frac{\partial }{\partial {r_i}^{(0)}} :=  (A_{T_0} r_i) \cdot  \nabla_{x^{(0)}} = \sum_{j_0=1}^d (A_{T_0} r_i)_{j_0} \frac{\partial}{\partial x_{j_0}^{(0)}}, \quad i \in \{ 1 : d \},
\end{align*}
where ${A}_{T_0} \in O(d)$ is the orthogonal matrix defined in Eq.~\eqref{affine10}. For a multi-index $\beta = (\beta_1,\ldots,\beta_d) \in \mathbb{N}_0^d$, we use the notation
\begin{align*}
\displaystyle
\partial^{\beta}_{r^{(0)}} := \frac{\partial^{|\beta|}}{ (\partial r_1^{(0)})^{\beta_1} \ldots (\partial r_d^{(0)})^{\beta_d}}.
\end{align*}
\end{defi}

%\begin{Rem}
%The vectors $r_i$, $i \in \{ 1,\ldots,d \}$ defined in Definition \ref{defi332} are unit vectors. 
%\end{Rem}

\begin{Note} \label{defi=theta}
Recall that
\begin{align*}
\displaystyle
&|s| \leq 1, \ \alpha_2 \leq \alpha_1 \quad \text{if $d=2$},\\
&|s_1|\leq 1, \ |s_{21}| \leq 1, \  \alpha_2 \leq \alpha_3\leq \alpha_1 \quad \text{if $d=3$}.
\end{align*}
When $d=3$, if Assumption \ref{ass1} is imposed, there exists a positive constant $M$, independent of $h_T$, such that $|s_{22}| \leq M \frac{\alpha_2 t_1}{\alpha_3}$. Thus, if $d=2$, we have
\begin{align*}
\displaystyle
%&\alpha_1 | [\widetilde{A}]_{j1} | \leq (\mathcal{H}^{(2)})_j, \quad \alpha_2 | [\widetilde{A}]_{j2} | \leq (\mathcal{H}^{(2)} )_j, \quad j=1,2,
&\alpha_1 | [\widetilde{A}]_{j1} | \leq \mathscr{H}_j, \quad \alpha_2 | [\widetilde{A}]_{j2} | \leq \mathscr{H}_j, \quad j=1,2,
%&\alpha_i | [\widetilde{A}]_{ji} | \leq (\mathcal{H}^{(2)}_i)_j, \quad i,j =1,2,
\end{align*}
and if $d=3$, for $\widetilde{A} \in \{ \widetilde{A}_1 , \widetilde{A}_2  \}$ and $j=1,2,3$, we have
\begin{align*}
\displaystyle
%&\alpha_1 | [\widetilde{A}]_{j1} | \leq ( \mathcal{H}^{(3)} )_j, \quad \alpha_2 | [\widetilde{A}]_{j2} | \leq (\mathcal{H}^{(3)})_j, \\
%&\alpha_3 | [\widetilde{A}]_{j3} | \leq \max \{ 1,M\} (\mathcal{H}^{(3)})_j,  \quad j=1,2,3.
&\alpha_1 | [\widetilde{A}]_{j1} | \leq \mathscr{H}_j, \quad \alpha_2 | [\widetilde{A}]_{j2} | \leq \mathscr{H}_j, \quad \alpha_3 | [\widetilde{A}]_{j3} | \leq \max \{ 1,M\} \mathscr{H}_j,  \quad j=1,2,3.
%&\alpha_1 | [\widetilde{A}]_{j1} | \leq  (\mathcal{H}^{(3)}_1)_j, \quad \alpha_2 | [\widetilde{A}]_{j2} | \leq  (\mathcal{H}^{(3)}_2)_j, \quad \alpha_3 | [\widetilde{A}]_{j3} | \leq \max \{ 1,M\}  (\mathcal{H}^{(3)}_3)_j.
\end{align*}
\begin{comment}
Furthermore, we define a directional derivative as
\begin{align*}
\displaystyle
%\frac{\partial }{\partial \mathcal{H}_i} := (  \Abs{( A_{T_0 } )} \mathcal{H}^{(d)}_i ) \cdot \nabla_{x^{(0)}}  = \sum_{j_0=1}^d  (  \Abs{( A_{T_0 } )} \mathcal{H}^{(d)}_i )_{j_0}  \frac{\partial}{\partial x_{j_0}^{(0)}},
\frac{\partial }{\partial \mathcal{H}_i} := (  \Abs{( A_{T_0 } )} \mathcal{H}^{(d)}_i ) \cdot \nabla_{x^{(0)}}  = \sum_{j_0=1}^d  (  \Abs{( A_{T_0 } )} \mathcal{H}^{(d)}_i )_{j_0}  \frac{\partial}{\partial x_{j_0}^{(0)}},
\end{align*}
where we define "$ \Abs{( A_{T_0 } ) }$" as
\begin{align*}
\displaystyle
[ \Abs{( A_{T_0 } ) }]_{ij} := | [A_{T_0 }]_{ij} |, \quad i,j \in \{ 1:d \}.
\end{align*}
For a multi-index $\beta = (\beta_1,\ldots,\beta_d) \in \mathbb{N}_0^d$, we use the following notation:
\begin{align*}
\displaystyle
\partial^{\beta}_{\mathcal{H}} := \frac{\partial^{|\beta|}}{ (\partial \mathcal{H}_1)^{\beta_1} \ldots (\partial \mathcal{H}_d)^{\beta_d}}.
\end{align*}
\end{comment}

\end{Note}

\begin{Note}
We use the following calculations in Lemma \ref{int=lem5}. For any multi-indices $\beta$ and $\gamma$, we have
\begin{align*}
\displaystyle
\partial^{\beta + \gamma}_{\hat{x}} &= \frac{\partial^{|\beta| + |\gamma|}}{\partial \hat{x}_1^{\beta_1} \cdots \partial \hat{x}_d^{\beta_d} \partial \hat{x}_1^{\gamma_1} \cdots \partial \hat{x}_d^{\gamma_d}} \notag\\
&\hspace{-1.5cm}  =  \underbrace{\sum_{i_1^{(1)} = 1}^d \alpha_1 [\widetilde{A}]_{i_1^{(1)} 1} \cdots   \sum_{i_{\beta_1}^{(1)} = 1}^d \alpha_1 [\widetilde{A}]_{i_{\beta_1}^{(1)} 1}   }_{\beta_1 \text{times}} \cdots \underbrace{ \sum_{i_1^{(d)}  = 1}^d \alpha_d [\widetilde{A}]_{i_1^{(d)} d}  \cdots \sum_{i_{\beta_d}^{(d)} = 1}^d \alpha_d [\widetilde{A}]_{i_{\beta_d}^{(d)} d}  }_{\beta_d \text{times}}  \notag \\
&\hspace{-1.2cm}  \underbrace{\sum_{j_1^{(1)}  = 1}^d \alpha_1 [\widetilde{A}]_{j_1^{(1)} 1} \cdots  \sum_{j_{\gamma_1}^{(1)}= 1}^d \alpha_1 [\widetilde{A}]_{j_{\gamma_1}^{(1)} 1} }_{\gamma_1 \text{times}}  \cdots  \underbrace{\sum_{j_1^{(d)} = 1}^d \alpha_d [\widetilde{A}]_{j_1^{(d)} d}  \cdots \sum_{j_{\gamma_d}^{(d)} = 1}^d \alpha_d [\widetilde{A}]_{j_{\gamma_d}^{(d)} d}  }_{\gamma_d \text{times}} \notag \\
&\hspace{-1.2cm}  \underbrace{\frac{\partial^{\beta_1}}{\partial {x}_{i_1^{1)}} \cdots \partial {x}_{i_{\beta_1}^{(1)}}}}_{\beta_1 \text{times}} \cdots \underbrace{\frac{\partial^{\beta_d}}{\partial {x}_{i_1^{(d)}} \cdots \partial {x}_{i_{\beta_d}^{(d)}}}}_{\beta_d \text{times}} \underbrace{ \frac{\partial^{\gamma_1}}{ \partial {x}_{j_1^{(1)}} \cdots \partial {x}_{j_{\gamma_1}^{(1)}} } }_{\gamma_1 \text{times}} \cdots \underbrace{\frac{\partial^{\gamma_d}}{ \partial {x}_{j_1^{(d)}} \cdots \partial {x}_{j_{\gamma_d}^{(d)}} }}_{\gamma_d \text{times}}.
\end{align*}

Let $\hat{\varphi} \in \mathcal{C}^{\ell}(\widehat{T})$ with $\tilde{\varphi} = \hat{\varphi} \circ \widehat{\Phi}^{-1}$ and ${\varphi} = \tilde{\varphi} \circ \widetilde{\Phi}^{-1}$. Then, for $1 \leq i \leq d$, 
\begin{align*}
\displaystyle
\left| \frac{\partial \hat{\varphi}}{\partial \hat{x}_i  } \right|
&\leq  \sum_{i_1^{(1)}=1}^d  \alpha_i  \left| [\widetilde{A}]_{i_1^{(1)} i}  \right|  \left|\frac{\partial \varphi}{\partial x_{i_1^{(1)}}}  \right| 
\leq 
\begin{cases}	
 \alpha_i   \| \widetilde{A}\|_{\max}  \sum_{i_1^{(1)}=1}^d \left| \frac{\partial \varphi}{\partial x_{i_1^{(1)}}} \right|  \quad \text{or},\\
 c \sum_{i_1^{(1)}=1}^d  \mathscr{H}_{i_1^{(1)}}  \left|   \frac{\partial \varphi}{\partial x_{i_1^{(1)}}} \right|,
\end{cases}
\end{align*}
and for $1 \leq i,j \leq d$, 
\begin{align*}
\displaystyle
\left| \frac{\partial^2 \hat{\varphi}}{\partial \hat{x}_i \partial \hat{x}_j } \right|
&= \left | \sum_{i_1^{(1)}, j_1^{(1)}=1}^d  \alpha_i \alpha_j  [\widetilde{A}]_{i_1^{(1)} i} [\widetilde{A}]_{j_1^{(1)} j}  \frac{\partial^2 \varphi}{\partial x_{i_1^{(1)}} \partial x_{j_1^{(1)}}} \right | \\
&\leq
\begin{cases}	
\alpha_i \alpha_j \| \widetilde{A}\|_{\max}^2  \sum_{i_1^{(1)}, j_1^{(1)}=1}^d \left |  \frac{\partial^2 \varphi}{\partial x_{i_1^{(1)}} \partial x_{j_1^{(1)}}} \right | \quad \text{or}, \\
\alpha_j  \sum_{j_1^{(1)}=1}^d  |[\widetilde{A}]_{j_1^{(1)} j}|   \left|  \sum_{i_1^{(1)} = 1}^d  \alpha_i [\widetilde{A}]_{i_1^{(1)} i} \frac{\partial^2 \varphi}{\partial x_{i_1^{(1)}} \partial x_{j_1^{(1)}}}   \right|  \\
\quad \leq c \alpha_j \| \widetilde{A} \|_{\max} \sum_{j_1^{(1)}=1}^d \sum_{i_1^{(1)} = 1}^d  \mathscr{H}_{i_1^{(1)}} \left | \frac{\partial^2 \varphi}{\partial x_{i_1^{(1)}} \partial x_{j_1^{(1)}}}   \right| \quad \text{or}, \\
c \sum_{i_1^{(1)} = 1}^d \sum_{j_1^{(1)}=1}^d  \mathscr{H}_{i_1^{(1)}}  \mathscr{H}_{j_1^{(1)}} \left | \frac{\partial^2 \varphi}{\partial x_{i_1^{(1)}} \partial x_{j_1^{(1)}}}   \right|.
\end{cases}
\end{align*}

\end{Note}

\begin{lem} \label{int=lem5}
Suppose that Assumption \ref{ass1} is imposed. Let $m \in \mathbb{N}_0$, $\ell \in \mathbb{N}_0$ with $\ell \geq m$ and $p \in [0,\infty]$. Let  $\beta := (\beta_1,\ldots,\beta_d) \in \mathbb{N}_0^d$ and $\gamma := (\gamma_1,\ldots,\gamma_d) \in \mathbb{N}_0^d$ be multi-indices with $|\beta| = m$ and $|\gamma| = \ell - m$. Then, for any $\hat{\varphi} \in W^{\ell,p}(\widehat{T})$ with $\tilde{\varphi} = \hat{\varphi} \circ \widehat{\Phi}^{-1}$ and ${\varphi} = \tilde{\varphi} \circ \widetilde{\Phi}^{-1}$, it holds that
\begin{align}
\displaystyle
\| \partial^{\beta} \partial^{\gamma} \hat{\varphi} \|_{L^p(\widehat{T})}
&\leq C_2^{SA} |\det({A}_T)|^{-\frac{1}{p}} \| \widetilde{A} \|^m_2 \alpha^{\beta} \sum_{|\epsilon| =  |\gamma|} \mathscr{H}^{\varepsilon} | \partial^{\epsilon} \varphi |_{W^{m,p}(T)}, \label{sa29}
\end{align}
where $C_2^{SA}$ is a constant that is independent of $T$ and $\widetilde{T}$. Here, for $p = \infty$ and any positive real $x$, $x^{- \frac{1}{p}} = 1$.
\begin{comment}
When $p = \infty$, for $ \hat{\varphi} \in W^{\ell,\infty}(\widehat{T})$ with $\tilde{\varphi} = \hat{\varphi} \circ \widehat{\Phi}^{-1}$ and ${\varphi} = \tilde{\varphi} \circ \widetilde{\Phi}^{-1}$,  it holds that
\begin{align}
\displaystyle
\| \partial^{\beta} \partial^{\gamma} \hat{\varphi} \|_{L^{\infty}(\widehat{T})}
&\leq C_2^{SA,\infty} \| \widetilde{A} \|^m_2 \alpha^{\beta} \sum_{|\epsilon| =  |\gamma|} \mathscr{H}^{\varepsilon} | \partial^{\epsilon} \varphi |_{W^{m,\infty}(T)}, \label{sa34=infty}
\end{align}
where $C_2^{SA,\infty}$ is a constant that are independent of $T$ and $\widetilde{T}$. 
\end{comment}
\end{lem}

\begin{pf*}
Let $\varepsilon = (\varepsilon_1,\ldots,\varepsilon_d) \in \mathbb{N}_0^d$ and $\delta = (\delta_1,\ldots,\delta_d) \in \mathbb{N}_0^d$ be multi-indies with $|\varepsilon| = |\gamma|$ and $|\delta| = |\beta|$. Let $p \in [1,\infty)$. Because the space $\mathcal{C}^{\ell}(\widehat{T})$ is dense in the space ${W}^{\ell,p}(\widehat{T})$, we show that Eq.~\eqref{sa29} holds for $\hat{\varphi} \in \mathcal{C}^{\ell}(\widehat{T})$ with $\tilde{\varphi} = \hat{\varphi} \circ \widehat{\Phi}^{-1}$ and ${\varphi} = \tilde{\varphi} \circ \widetilde{\Phi}^{-1}$. Through a simple calculation, we obtain
\begin{align}
\displaystyle
|\partial^{\beta + \gamma} \hat{\varphi} | &= \left| \frac{\partial^{\ell} \hat{ \varphi}}{\partial \hat{x}_1^{\beta_1} \cdots \partial \hat{x}_d^{\beta_d} \partial \hat{x}_1^{\gamma_1} \cdots \partial \hat{x}_d^{\gamma_d}} \right| \notag\\
&\hspace{-1.5cm}  \leq c \alpha^{\beta}  \| \widetilde{A} \|_{\max}^{|\beta|}  \underbrace{\sum_{i_1^{(1)} = 1}^d  \cdots   \sum_{i_{\beta_1}^{(1)} = 1}^d }_{\beta_1 \text{times}} \cdots \underbrace{ \sum_{i_1^{(d)} = 1}^d  \cdots \sum_{i_{\beta_d}^{(d)} = 1}^d  }_{\beta_d \text{times}}  \underbrace{\sum_{j_1^{(1)} = 1}^d  \cdots  \sum_{j_{\gamma_1}^{(1)} = 1}^d  }_{\gamma_1 \text{times}}  \cdots   \underbrace{\sum_{j_1^{(d)} = 1}^d \cdots \sum_{j_{\gamma_d}^{(d)} = 1}^d }_{\gamma_d \text{times}} \notag \\
&\hspace{-1.2cm}  \underbrace{\mathscr{H}_{j_1^{(1)}} \cdots \mathscr{H}_{j_{\varepsilon_1}^{(1)}}}_{\gamma_1 \text{times}} \cdots  \underbrace{ \mathscr{H}_{j_1^{(d)}} \cdots \mathscr{H}_{j_{\varepsilon_d}^{(d)}}}_{\gamma_d \text{times}} \notag \\
&\hspace{-1.2cm} \Biggl | \underbrace{\frac{\partial^{\beta_1}}{\partial {x}_{i_1^{1)}} \cdots \partial {x}_{i_{\beta_1}^{(1)}}}}_{\beta_1 \text{times}} \cdots \underbrace{\frac{\partial^{\beta_d}}{\partial {x}_{i_1^{(d)}} \cdots \partial {x}_{i_{\beta_d}^{(d)}}}}_{\beta_d \text{times}} \underbrace{ \frac{\partial^{\gamma_1}}{ \partial {x}_{j_1^{(1)}} \cdots \partial {x}_{j_{\gamma_1}^{(1)}} } }_{\gamma_1 \text{times}} \cdots \underbrace{\frac{\partial^{\gamma_d}}{ \partial {x}_{j_1^{(d)}} \cdots \partial {x}_{j_{\gamma_d}^{(d)}} }}_{\gamma_d \text{times}} \varphi \Biggr |  \notag \\
&\hspace{-1.5cm}  \leq c \alpha^{\beta}  \| \widetilde{A} \|_{\max}^{|\beta|}   \sum_{|\delta| = |\beta|} \sum_{|\varepsilon| = |\gamma|} \mathscr{H}^{\varepsilon}  | { \partial^{\delta} \partial^{\varepsilon} \varphi}  |. \notag
\end{align}
Using Eq.~\eqref{matrix22}, we then have that
\begin{align*}
\displaystyle
\int_{\widehat{T}} | \partial^{\beta} \partial^{ \gamma} \hat{\varphi}|^p d \hat{x}
&\leq c  \| \widetilde{A}  \|_2^{mp} \alpha^{\beta p}  \sum_{|\delta| = |\beta|} \sum_{|\epsilon| = |\gamma|} \mathscr{H}^{\varepsilon p} \int_{\widehat{T}} | { \partial^{\delta} \partial^{\varepsilon} \varphi}  |^p d \hat{x} \\
&= c |\det(A_T)|^{-1} \| \widetilde{A}  \|_2^{mp} \alpha^{\beta p}  \sum_{|\delta| = |\beta|} \sum_{|\epsilon| = |\gamma|} \mathscr{H}^{\varepsilon p} \int_{{T}} |{\partial^{\delta} \partial^{\varepsilon} \varphi}  |^p d {x}.
\end{align*}
Therefore, using \eqref{jensen}, we obtain
\begin{align*}
\displaystyle
\|  \partial^{\beta} \partial^{\gamma} \hat{\varphi}\|_{L^p(\widehat{T})}
&\leq c |\det(A_T)|^{-\frac{1}{p}} \| \widetilde{A}  \|_2^{m} \alpha^{\beta} \sum_{|\epsilon| =  |\gamma|} \mathscr{H}^{\varepsilon} | \partial^{\epsilon} \varphi |_{W^{m,p}(T)},
\end{align*}
which recovers Eq.~\eqref{sa29}. 

We consider the case in which $p = \infty$. A function $\varphi \in W^{\ell,\infty}(T)$ belongs to the space $W^{\ell,p}(T)$ for any $p \in [1,\infty)$. Therefore, it holds that $\hat{\varphi} \in W^{\ell,p}(\widehat{T})$ for any $p \in [1,\infty)$, and thus,
\begin{align}
\displaystyle
\|  \partial^{\beta} \partial^{\gamma} \hat{\varphi}\|_{L^p(\widehat{T})}
&\leq c |\det(A_T)|^{-\frac{1}{p}} \| \widetilde{A}  \|_2^{m} \alpha^{\beta} \sum_{|\epsilon| =  |\gamma|} \mathscr{H}^{\varepsilon} | \partial^{\epsilon} \varphi |_{W^{m,p}(T)} \notag \\
&\leq c \| \widetilde{A}  \|_2^{m} \alpha^{\beta} \sum_{|\epsilon| =  |\gamma|} \mathscr{H}^{\varepsilon} | \partial^{\epsilon} \varphi |_{W^{m,\infty}(T)} \< \infty. \label{infty32}
\end{align}
This implies that the function $\partial^{\beta} \partial^{\gamma} \hat{\varphi}$ is in the space $L^{\infty}(\widehat{T})$. Inequality  \eqref{sa29} for $p=\infty$ is obtained by taking the limit $p \to \infty$ in Eq.~\eqref{infty32} on the basis that $\lim_{p \to \infty} \| \cdot \|_{L^p(\widehat{T})} = \| \cdot\|_{L^{\infty}(\widehat{T})}$.
\qed
\end{pf*}

\begin{Rem} \label{ex=01}
In inequality \eqref{sa29}, it is possible to obtain the estimates in $T_0$ by specifically determining the matrix $A_{T_0}$.

Let $\ell=2$, $m=1$, and $p=q=2$. Recall that
\begin{align*}
\displaystyle
\Phi_{T_0}: T \ni x \mapsto x^{(0)} := {A}_{T_0} x + b_{T_0} \in T_0.
\end{align*}
For ${\varphi} \in \mathcal{C}^{2}({T})$ with $\varphi_0 = \varphi \circ \Phi_{T_0}^{-1}$ and $1 \leq i,j \leq d$, 
we have
\begin{align*}
\displaystyle
\left| \frac{\partial^2 {\varphi}}{\partial {x}_i  {\partial {x}_j}} ({x}) \right|
&= \left| \sum_{i_1^{(1)} , j_1^{(1)}=1}^2  [A_{T_0}]_{i_1^{(1)} i} [A_{T_0}]_{j_1^{(1)} j} \frac{\partial^2 \varphi_0}{\partial x_{i_1^{(1)}}^{(0)} \partial x_{j_1^{(1)}}^{(0)}} (x)  \right|.
\end{align*}

Let $d=2$. We define the matrix $A_{T_0}$ as 
\begin{align*}
\displaystyle
A_{T_0} := 
\begin{pmatrix}
\cos \frac{\pi}{2}  & - \sin \frac{\pi}{2} \\
 \sin \frac{\pi}{2} & \cos \frac{\pi}{2}
\end{pmatrix}.
\end{align*}
Because $ \| A_{T_0} \|_{\max} = 1$, we have
\begin{align*}
\displaystyle
\left| \frac{\partial^2 {\varphi}}{\partial {x}_i  {\partial {x}_j}} ({x}) \right|
&\leq \left|  \frac{\partial^2 \varphi_0}{\partial x_{i+1}^{(0)} \partial x_{j+1}^{(0)}} (x)  \right|,
\end{align*}
where the indices $i$, $i+1$ and $j$, $j+1$ must be evaluated modulo 2.  Because $|\det (A_{T_0})| = 1$, it holds that
\begin{align*}
\displaystyle
\left \| \frac{\partial^2 {\varphi}}{\partial {x}_i  {\partial {x}_j}} \right \|_{L^2(T)} \leq  \left \|  \frac{\partial^2 \varphi_0}{\partial x_{i_{i+1}}^{(0)} \partial x_{j+1}^{(0)}}  \right \|_{L^2(T_0)}.
\end{align*}
We then have
\begin{align*}
\displaystyle
\sum_{j = 1}^2 \mathscr{H}_j \left | \frac{\partial {\varphi}}{ {\partial {x}_j}} \right |_{H^1(T)}
&\leq \sum_{j = 1}^2 \mathscr{H}_j \left |  \frac{\partial \varphi_0}{ \partial x_{j+1}^{(0)}}  \right |_{H^1(T_0)},
\end{align*}
where the indices $j$, $j+1$ must be evaluated modulo 2. 

We define the matrix $A_{T_0}$ as 
\begin{align*}
\displaystyle
A_{T_0} := 
\begin{pmatrix}
\cos \frac{\pi}{4}  & - \sin \frac{\pi}{4} \\
 \sin \frac{\pi}{4} & \cos \frac{\pi}{4}
\end{pmatrix}.
\end{align*}
We then have
\begin{align*}
\displaystyle
\left| \frac{\partial^2 {\varphi}}{\partial {x}_i  {\partial {x}_j}} ({x}) \right|
&\leq \frac{1}{\sqrt{2}}  \sum_{i_1^{(1)} , j_1^{(1)}=1}^2  \left|  \frac{\partial^2 \varphi_0}{\partial x_{i_1^{(1)}}^{(0)} \partial x_{j_1^{(1)}}^{(0)}} (x)  \right|,
\end{align*}
which leads to
\begin{align*}
\displaystyle
\left \| \frac{\partial^2 {\varphi}}{\partial {x}_i  {\partial {x}_j}} \right \|_{L^2(T)}^2 \leq  c \sum_{i_1^{(1)} , j_1^{(1)}=1}^2 \left \|  \frac{\partial^2 \varphi_0}{\partial x_{i_1^{(1)}}^{(0)} \partial x_{j_1^{(1)}}^{(0)}}  \right \|_{L^2(T_0)}^2 \leq c |\varphi_0|^2_{H^2(T_0)}.
\end{align*}
Using \eqref{jensen}, we than have that
\begin{align*}
\displaystyle
\sum_{j = 1}^2 \mathscr{H}_j \left | \frac{\partial {\varphi}}{ {\partial {x}_j}} \right |_{H^1(T)}
&\leq \sum_{j = 1}^2 \mathscr{H}_j |\varphi_0|_{H^2(T_0)} \leq c h_{T_0} |\varphi_0|_{H^2(T_0)}.
\end{align*}
In this case, anisotropic interpolation error estimates cannot be obtained.

\end{Rem}

\begin{Note}
We use the following calculations in Lemma \ref{int=lem6}. For any multi-indices $\beta$ and $\gamma$, we have
\begin{align*}
\displaystyle
\partial^{\beta + \gamma}_{\hat{x}} &= \frac{\partial^{|\beta| + |\gamma|}}{\partial \hat{x}_1^{\beta_1} \cdots \partial \hat{x}_d^{\beta_d} \partial \hat{x}_1^{\gamma_1} \cdots \partial \hat{x}_d^{\gamma_d}} \notag\\
&\hspace{-1.5cm}  =  \underbrace{\sum_{i_1^{(1)},i_1^{(0,1)} = 1}^d \alpha_1 [\widetilde{A}]_{i_1^{(1)} 1} [A_{T_0}]_{i_1^{(0,1)} i_1^{(1)}} \cdots   \sum_{i_{\beta_1}^{(1)},i_{\beta_1}^{(0,1)} = 1}^d \alpha_1 [\widetilde{A}]_{i_{\beta_1}^{(1)} 1} [A_{T_0}]_{i_{\beta_1}^{(0,1)} i_{\beta_1}^{(1)}}   }_{\beta_1 \text{times}} \cdots \notag \\
&\hspace{-1.2cm} \underbrace{ \sum_{i_1^{(d)} , i_1^{(0,d)} = 1}^d \alpha_d [\widetilde{A}]_{i_1^{(d)} d} [A_{T_0}]_{i_1^{(0,d)} i_1^{(d)}}  \cdots \sum_{i_{\beta_d}^{(d)} , i_{\beta_d}^{(0,d)}= 1}^d \alpha_d [\widetilde{A}]_{i_{\beta_d}^{(d)} d} [A_{T_0}]_{i_{\beta_d}^{(0,d)} i_{\beta_d}^{(d)}}  }_{\beta_d \text{times}}  \notag \\
&\hspace{-1.2cm}  \underbrace{\sum_{j_1^{(1)} , j_1^{(0,1)} = 1}^d \alpha_1 [\widetilde{A}]_{j_1^{(1)} 1} [A_{T_0}]_{j_1^{(0,1)} j_1^{(1)}}  \cdots  \sum_{j_{\gamma_1}^{(1)} , j_{\gamma_1}^{(0,1)}= 1}^d \alpha_1 [\widetilde{A}]_{j_{\gamma_1}^{(1)} 1} [A_{T_0}]_{j_{\gamma_1}^{(0,1)} j_{\gamma_1}^{(1)}} }_{\gamma_1 \text{times}}  \cdots   \notag \\
&\hspace{-1.2cm}  \underbrace{\sum_{j_1^{(d)} , j_1^{(0,d)} = 1}^d \alpha_d [\widetilde{A}]_{j_1^{(d)} d} [A_{T_0}]_{j_1^{(0,d)} j_1^{(d)}}  \cdots \sum_{j_{\gamma_d}^{(d)} , j_{\gamma_d}^{(0,d)}= 1}^d \alpha_d [\widetilde{A}]_{j_{\gamma_d}^{(d)} d} [A_{T_0}]_{j_{\gamma_d}^{(0,d)} j_{\gamma_d}^{(d)}}  }_{\gamma_d \text{times}} \notag \\
&\hspace{-1.2cm}  \underbrace{\frac{\partial^{\beta_1}}{\partial {x}_{i_1^{(0,1)}}^{(0)} \cdots \partial {x}_{i_{\beta_1}^{(0,1)}}^{(0)}}}_{\beta_1 \text{times}} \cdots \underbrace{\frac{\partial^{\beta_d}}{\partial {x}_{i_1^{(0,d)}}^{(0)} \cdots \partial {x}_{i_{\beta_d}^{(0,d)}}^{(0)}}}_{\beta_d \text{times}}  \underbrace{ \frac{\partial^{\gamma_1}}{ \partial {x}_{j_1^{(0,1)}}^{(0)} \cdots \partial {x}_{j_{\gamma_1}^{(0,1)}}^{(0)} } }_{\gamma_1 \text{times}} \cdots \underbrace{\frac{\partial^{\gamma_d}}{ \partial {x}_{j_1^{(0,d)}}^{(0)} \cdots \partial {x}_{j_{\gamma_d}^{(0,d)}}^{(0)} }}_{\gamma_d \text{times}}.
\end{align*}

Let $\hat{\varphi} \in \mathcal{C}^{\ell}(\widehat{T})$ with $\tilde{\varphi} = \hat{\varphi} \circ \widehat{\Phi}^{-1}$, ${\varphi} = \tilde{\varphi} \circ \widetilde{\Phi}^{-1}$ and $\varphi_0 = \varphi \circ \Phi_{T_0}^{-1}$. Then, for $1 \leq i \leq d$, 
\begin{align*}
\displaystyle
\left| \frac{\partial \hat{\varphi}}{\partial \hat{x}_i  } \right|
&= \left|  \sum_{i_1^{(1)}=1}^d  \sum_{i_1^{(0,1)} =1}^d \alpha_i  [\widetilde{A}]_{i_1^{(1)} i} [A_{T_0}]_{i_1^{(0,1)} i_1^{(1)}}  \frac{\partial \varphi_0}{\partial x_{i_1^{(0,1)}}^{(0)}}  \right| \\
&= \alpha_i   \left |  \sum_{i_1^{(1)}=1}^d  \sum_{i_1^{(0,1)}=1}^d  [A_{T_0}]_{i_1^{(0,1)} i_1^{(1)}} (r_i)_{i_1^{(1)}}  \frac{\partial \varphi_0}{\partial x_{i_1^{(0,1)}}^{(0)}} \right| = \alpha_i \left |  \frac{\partial \varphi_0}{\partial r_{i}^{(0)}} \right| \\
&\leq \alpha_i \| \widetilde{A} \|_{\max} \| A_{T_0} \|_{\max}  \sum_{i_1^{(1)}=1}^d  \sum_{i_1^{(0,1)}=1}^d \left | \frac{\partial \varphi_0}{\partial x_{i_1^{(0,1)}}^{(0)}} \right|,
\end{align*}
and for $1 \leq i,j \leq d$, 
\begin{align*}
\displaystyle
\left| \frac{\partial^2 \hat{\varphi}}{\partial \hat{x}_i \partial \hat{x}_j } \right|
&= \Biggl | \sum_{i_1^{(1)}, j_1^{(1)}=1}^d  \sum_{i_1^{(0,1)}, j_1^{(0,1)}=1}^d \alpha_i  \alpha_j  [\widetilde{A}]_{i_1^{(1)} i} [\widetilde{A}]_{j_1^{(1)} j}  \\
&\quad \quad [A_{T_0}]_{i_1^{(0,1)} i_1^{(1)}} [A_{T_0}]_{j_1^{(0,1)} j_1^{(1)}}  \frac{\partial^2 \varphi_0}{\partial x_{i_1^{(0,1)}}^{(0)} \partial x_{j_1^{(0,1)}}^{(0)}} \Biggr | = \alpha_i  \alpha_j \left |  \frac{\partial^2 \varphi_0}{\partial r_{i}^{(0)} \partial r_{j}^{(0)}} \right| \\
&\leq \alpha_i  \alpha_j \sum_{j_1^{(1)} =1}^d | [\widetilde{A}]_{j_1^{(1)} j} |    \Biggl | \sum_{j_1^{(0,1)}=1}^d [A_{T_0}]_{j_1^{(0,1)} j_1^{(1)}}  \frac{\partial^2 \varphi_0}{\partial r_{i}^{(0)} \partial x_{j_1^{(0,1)}}^{(0)}} \Biggr | \\
&\leq \alpha_i  \alpha_j  \| \widetilde{A} \|_{\max}   \| {A}_{T_0} \|_{\max} \sum_{j_1^{(0,1)}=1}^d \Biggl | \frac{\partial^2 \varphi_0}{\partial r_{i}^{(0)} \partial x_{j_1^{(0,1)}}^{(0)}} \Biggr | \\
&\leq \alpha_i  \alpha_j  \| \widetilde{A} \|_{\max}^2   \| {A}_{T_0} \|_{\max}^2  \sum_{i_1^{(0,1)}, j_1^{(0,1)}=1}^d \left| \frac{\partial^2 \varphi_0}{\partial x_{i_1^{(0,1)}}^{(0)} \partial x_{j_1^{(0,1)}}^{(0)}} \right|.
\end{align*}
\end{Note}

If Assumption \ref{ass1} is not imposed, the estimates corresponding to Lemma \ref{int=lem5} are as follows.

\begin{lem} \label{int=lem6}
Let $m \in \mathbb{N}_0$, $\ell \in \mathbb{N}_0$ with $\ell \geq m$, and $p \in [0,\infty]$. Let  $\beta := (\beta_1,\ldots,\beta_d) \in \mathbb{N}_0^d$ and $\gamma := (\gamma_1,\ldots,\gamma_d) \in \mathbb{N}_0^d$ be multi-indices with $|\beta| = m$ and $|\gamma| = \ell - m$. Then, for any $\hat{\varphi} \in W^{\ell,p}(\widehat{T})$ with $\tilde{\varphi} = \hat{\varphi} \circ \widehat{\Phi}^{-1}$, ${\varphi} = \tilde{\varphi} \circ \widetilde{\Phi}^{-1}$, and $\varphi_0 = \varphi \circ \Phi_{T_0}^{-1}$, it holds that
\begin{align}
\displaystyle
\| \partial^{\beta} \partial^{\gamma}  \hat{\varphi}\|_{L^p(\widehat{T})}
&\leq  C_3^{SA} |\det({A}_T)|^{-\frac{1}{p}} \| \widetilde{{A}}  \|_2^{m} \alpha^{\beta} \sum_{|\epsilon| =  |\gamma|} {\alpha}^{\varepsilon} | \partial_{r^{(0)}}^{\epsilon} \varphi_0 |_{W^{m,p}(T_0)}, \label{sa32}
\end{align}
where $C_3^{SA}$ is a constant that is independent of $T_0$ and $\widetilde{T}$. Here, for $p = \infty$ and any positive real $x$, $x^{- \frac{1}{p}} = 1$.

\begin{comment}
When $p = \infty$, for $ \hat{\varphi} \in W^{\ell,\infty}(\widehat{T})$ with $\tilde{\varphi} = \hat{\varphi} \circ \widehat{\Phi}^{-1}$, ${\varphi} = \tilde{\varphi} \circ \widetilde{\Phi}^{-1}$ and $\varphi_0 = \varphi \circ \Phi_{T_0}^{-1}$,  it holds that
\begin{align}
\displaystyle
\| \partial^{\beta} \partial^{\gamma}  \hat{\varphi}\|_{L^{\infty}(\widehat{T})}
&\leq  C_3^{SA,\infty} \| \widetilde{{A}}  \|_2^{m} \alpha^{\beta} \sum_{|\epsilon| =  |\gamma|} {\alpha}^{\varepsilon} | \partial_{r^{(0)}}^{\epsilon} \varphi_0 |_{W^{m,\infty}(T_0)}, \label{sa37=infty}
\end{align}
where $C_3^{SA,\infty}$ is a constant that are independent of $T_0$ and $\widetilde{T}$. 
\end{comment}
\end{lem}

\begin{pf*}
We follow the proof of Lemma \ref{int=lem5}. Let $p \in [1,\infty)$. Because the space $\mathcal{C}^{\ell}(\widehat{T})$ is dense in the space ${W}^{\ell,p}(\widehat{T})$, we show that Eq.~\eqref{sa32} holds for $\hat{\varphi} \in \mathcal{C}^{\ell}(\widehat{T})$ with $\tilde{\varphi} = \hat{\varphi} \circ \widehat{\Phi}^{-1}$, ${\varphi} = \tilde{\varphi} \circ \widetilde{\Phi}^{-1}$, and $\varphi_0 = \varphi \circ \Phi_{T_0}^{-1}$. For $1 \leq i,k \leq d$, 
\begin{align*}
\displaystyle
\left| \partial^{\beta + \gamma} \hat{\varphi} \right| 
&\leq c \alpha^{\beta}  \| \widetilde{A} \|_{\max}^{|\beta|}  \| {A}_{T_0} \|_{\max}^{|\beta|} \sum_{|\delta| = |\beta|}  \sum_{|\varepsilon| = |\gamma|} \alpha^{\varepsilon}  \left|  \partial^{\delta}  \partial_{r^{(0)}}^{\varepsilon} \varphi_0  \right|.
\end{align*}
Using Eqs.~\eqref{CN331c} and \eqref{matrix22}, we obtain Eq.~\eqref{sa32} for $p \in [1,\infty]$ by an argument analogous to that used for Lemma \ref{int=lem5}.
\qed
\end{pf*}

\section{Remarks on anisotropic interpolation analysis} \label{sec_mis}
We use the following Bramble--Hilbert-type lemma on anisotropic meshes proposed in \cite[Lemma 2.1]{Ape99}.

\begin{lem} \label{BH=lem7}
Let $D \subset \mathbb{R}^d$ with $d \in \{ 2,3\}$ be a connected open set that is star-shaped with respect to a ball $B$. Let $\gamma$ be a multi-index with $m := |\gamma|$ and $\varphi \in L^1(D)$ be a function with $\partial^{\gamma} \varphi \in W^{\ell -m,p}(D)$, where $\ell \in \mathbb{N}$, $m \in \mathbb{N}_0$, $0 \leq m \leq \ell$, and $p \in [1,\infty]$. Then, it holds that
\begin{align}
\displaystyle
\| \partial^{\gamma} (\varphi - Q^{(\ell)} \varphi) \|_{W^{\ell -m,p}(D)} \leq C^{BH} |\partial^{\gamma} \varphi|_{W^{\ell-m,p}(D)},  \label{lem7=1}
\end{align}
where $C^{BH}$ depends only on $d$, $\ell$, $\diam D$, and $\diam B$, and $ Q^{(\ell)} \varphi$ is defined as
\begin{align}
\displaystyle
(Q^{(\ell)} \varphi)(x) := \sum_{|\delta| \leq \ell -1} \int_B \eta(y) (\partial^{\delta}\varphi)(y) \frac{(x-y)^{\delta}}{\delta !} dy \in \mathcal{P}^{\ell -1},  \label{lem7=2}
\end{align}
where $\eta \in \mathcal{C}_0^{\infty}(B)$ is a given function with $\int_B \eta dx = 1$.
\end{lem}

As explained in the Introduction, there exist some mistakes in the proof of Theorem 2 of \cite{IshKobTsu}, and the statement is not valid in its original form. To clarify the following description, we explain the errors in the proof. Let $\widehat{T} \subset \mathbb{R}^2$ be the reference element defined in Section \ref{reference2d}. We set $k = m  = 1$, $\ell = 2$, and $p=2$. For $\hat{\varphi} \in H^2(\widehat{T})$, we set $\tilde{\varphi} := \hat{\varphi} \circ \widehat{\Phi}^{-1}$ and ${\varphi} := \tilde{\varphi} \circ \widetilde{\Phi}^{-1}$. Inequality \eqref{sa25} yields
\begin{align}
\displaystyle
|{\varphi} - I_T \varphi|_{H^{1}({T})}
&\leq c |\det(A_T)|^{\frac{1}{2}} \| \widetilde{A}^{-1} \|_2 \left( \sum_{i=1}^2  \alpha^{- 2}_i \|  \partial_{\hat{x}_i} ( \hat{\varphi} - I_{\widehat{T}} \hat{\varphi} ) \|^2_{L^2(\widehat{T})} \right)^{\frac{1}{2}}. \label{rem36}
\end{align}
The coefficient $\alpha_i^{-2}$ appears on the right-hand side of Eq.~\eqref{rem36}. In \cite[Theorem 2]{IshKobTsu}, we wrongly claimed that $\alpha_i^{-2}$ could be canceled out. In fact, a further assumption is required for this. Using  Eq.~\eqref{lem7=2} and the triangle inequality, we have
\begin{align*}
\displaystyle
 \|  \partial_{\hat{x}_i} ( \hat{\varphi} - I_{\widehat{T}} \hat{\varphi} ) \|^2_{L^2(\widehat{T})}
 &\leq 2 \|  \partial_{\hat{x}_i} ( \hat{\varphi} - Q^{(2)} \hat{\varphi} ) \|^2_{L^2(\widehat{T})} + 2 \|  \partial_{\hat{x}_i} ( Q^{(2)} \hat{\varphi}  -   I_{\widehat{T}} \hat{\varphi} ) \|^2_{L^2(\widehat{T})}.
\end{align*}
We use inequality \eqref{lem7=1} to obtain the target inequality \cite[Theorem 2]{IshKobTsu}. To this end, we have to show that
\begin{align}
\displaystyle
 \|  \partial_{\hat{x}_i} ( Q^{(2)} \hat{\varphi}  -   I_{\widehat{T}} \hat{\varphi} ) \|_{L^2(\widehat{T})}
 &\leq c \|  \partial_{\hat{x}_i} ( \hat{\varphi} - Q^{(2)} \hat{\varphi} ) \|_{H^1(\widehat{T})}. \label{rem37}
\end{align}
However, this is unlikely to hold because Eqs.~\eqref{int1} and \eqref{int513} yield
\begin{align*}
\displaystyle
 \|  \partial_{\hat{x}_i} ( Q^{(2)} \hat{\varphi}  -   I_{\widehat{T}} \hat{\varphi} ) \|_{L^2(\widehat{T})}
 &=  \|  \partial_{\hat{x}_i} (    I_{\widehat{T}} ( Q^{(2)} \hat{\varphi} )  -   I_{\widehat{T}} \hat{\varphi} ) \|_{L^2(\widehat{T})} \\
 &\leq c \| Q^{(2)} \hat{\varphi} -   \hat{\varphi} \|_{H^2(\widehat{T})} \leq c |\hat{\varphi}|_{H^2(\widehat{T})}.
\end{align*}
Using the classical scaling argument (see \cite[Lemma 1.101]{ErnGue04}), we have
\begin{align*}
\displaystyle
|\hat{\varphi}|_{H^2(\widehat{T})} \leq c |\det(A_T)|^{- \frac{1}{2}} \| \widetilde{A} \|_2 |\varphi|_{H^2(T)},
\end{align*}
which does not include the quantity $\alpha_i$. Therefore, the quantity $\alpha_i^{-1}$ in Eq.~\eqref{rem36} remains. Thus, the proof of \cite[Theorem 2]{IshKobTsu} is incorrect.

To overcome this problem, we use some results from previous studies \cite{ApeDob92,Ape99}. That is, we assume that there exists a linear functional $\mathscr{F}_1$ such that
\begin{align*}
\displaystyle
&\mathscr{F}_1 \in H^{1} (\widehat{T})^{\prime}, \\
&\mathscr{F}_1 ( \partial_{\hat{x}_i} (\hat{\varphi} - I_{\widehat{T}} \hat{\varphi}) ) = 0 \quad i = 1,2, \quad \forall \hat{\varphi} \in \mathcal{C}({\widehat{T}}): \ \partial_{\hat{x}_i} \hat{\varphi} \in H^{1} (\widehat{T}), \\
&\hat{\eta} \in \mathcal{P}^{1}, \quad \mathscr{F}_1(\partial_{\hat{x}_i} \hat{\eta}) = 0 \quad i = 1,2, \quad \Rightarrow \quad \partial_{\hat{x}_i} \hat{\eta} = 0.
\end{align*}
Because the polynomial spaces are finite-dimensional, all norms are equivalent; i.e., because $ | \mathscr{F}_1( \partial_{\hat{x}_i} ( \hat{\eta}  -   I_{\widehat{T}} \hat{\varphi} ) )|$ ($i=1,2$) is a norm on $\mathcal{P}^{0}$, we have that, for $i=1,2$,
\begin{align*}
\displaystyle
\| \partial_{\hat{x}_i} (\hat{\eta} - I_{\widehat{T}} \hat{\varphi}) \|_{L^2(\widehat{T})}
&\leq c  | \mathscr{F}_1 (\partial_{\hat{x}_i} ( \hat{\eta} - I_{\widehat{T}} \hat{\varphi} ) )| 
= c |\mathscr{F}_1 (\partial_{\hat{x}_i} ( \hat{\eta} - \hat{\varphi} ) )| \\
&\leq c \| \partial_{\hat{x}_i} ( \hat{\eta} - \hat{\varphi} ) \|_{H^1(\widehat{T})}.
\end{align*}
Setting $\hat{\eta} := Q^{(2)} \hat{\varphi}$, we obtain Eq.~\eqref{rem37}. Using inequality \eqref{lem7=1} yields
\begin{align*}
\displaystyle
 \|  \partial_{\hat{x}_i} ( \hat{\varphi} - I_{\widehat{T}} \hat{\varphi} ) \|^2_{L^2(\widehat{T})}
 &\leq c |  \partial_{\hat{x}_i} \hat{\varphi} |^2_{H^1(\widehat{T})},
\end{align*}
and so inequality \eqref{rem36} together with Eq.~\eqref{jensen} can be written as
\begin{align}
\displaystyle
|{\varphi} - I_T \varphi|_{H^{1}({T})}
&\leq c |\det(A_T)|^{\frac{1}{2}} \| \widetilde{A}^{-1} \|_2 \sum_{i,j=1}^2  \alpha^{- 1}_i \|  \partial_{\hat{x}_i}  \partial_{\hat{x}_j}  \hat{\varphi} \|_{L^2(\widehat{T})}. \label{rem38}
\end{align}
Inequality \eqref{sa29} yields
\begin{align}
\displaystyle
\|  \partial_{\hat{x}_i}  \partial_{\hat{x}_j} \hat{\varphi} \|_{L^2(\widehat{T})}
&\leq c |\det({A}_T)|^{-\frac{1}{2}} \| \widetilde{A} \|_2 \alpha_i \sum_{n=1}^2 \mathscr{H}_n \left|  \frac{ \partial \varphi}{\partial x_n} \right|_{H^{1}(T)}. \label{rem39}
\end{align}
Therefore, the quantity $\alpha_i^{-1}$ in Eq.~\eqref{rem38} and the quantity $\alpha_i$ in Eq.~\eqref{rem39} cancel out.

\section{Classical interpolation error estimates} \label{sec_cla}
The following embedding results hold.
\begin{thrso*}
 Let $d \geq 2$, $s \> 0$, and $p \in [1,\infty]$. Let $D \subset \mathbb{R}^d$ be a bounded open subset of $\mathbb{R}^d$. If $D$ is a Lipschitz set, we have that
\begin{align}
\displaystyle
W^{s,p}(D) \hookrightarrow
\begin{cases}
L^q(D) \quad \text{$\forall q \in [p , \frac{pd}{d - sp}]$ if $s p \< d$}, \\
L^q(D) \quad \text{$\forall q \in [p,\infty)$, if $sp=d$},\\
L^{\infty}(D) \cap \mathcal{C}^{0,\xi}(\overline{D}) \quad \text{$\xi = 1 - \frac{d}{s p}$ if $sp \> d$}.
\end{cases} \label{emmed}
\end{align} 
Furthermore, 
\begin{align}
\displaystyle
W^{s,p}(D) \hookrightarrow L^{\infty}(D) \cap \mathcal{C}^{0}(\overline{D}) \quad \text{(case $s=d$ and $p=1$)}. \label{emmed2}
\end{align} 
\end{thrso*}

\begin{pf*}
See, for example, \cite[Corollary B.43,  Theorem B.40]{ErnGue04}, \cite[Theorem 2.31]{ErnGue21a}, and the references therein.
\qed	
\end{pf*}

\begin{Rem}
 Let $s \> 0$ and $p \in [1,\infty]$ be such that
\begin{align*}
\displaystyle
s \> \frac{d}{p} \quad \text{if $p  \> 1$}, \quad s \geq d \quad \text{if $p=1$}.
\end{align*}
Then, it holds that  $W^{s,p}(D) \hookrightarrow \mathcal{C}^0(\overline{D})$.
\end{Rem}

 Using the new geometric parameter $H_{T_0}$, it is possible to deduce the classical interpolation error estimates; e.g., see \cite[Theorem 1.103]{ErnGue04} and \cite[Theorem 11.13]{ErnGue21a}.

\begin{th21*}
Let $1 \leq p \leq \infty$ and assume that there exists a nonnegative integer $k$ such that
\begin{align*}
\displaystyle
\mathcal{P}^{k} \subset \widehat{{P}} \subset W^{k+1,p}(\widehat{T}) \subset V(\widehat{T}).
\end{align*}
Let $\ell$ ($0 \leq \ell \leq k$) be such that $W^{\ell+1,p}(\widehat{T}) \subset V(\widehat{T})$ with the continuous embedding. Furthermore, assume that $\ell, m \in \mathbb{N} \cup \{ 0 \}$ and $p , q \in [1,\infty]$ such that $0 \leq m \leq \ell + 1$ and 
\begin{align}
\displaystyle
W^{\ell +1,p}(\widehat{T}) \hookrightarrow W^{m,q} (\widehat{T}). \label{cla40}
\end{align}
Then, for any $\varphi_0 \in W^{\ell+1,p}(T_0)$, it holds that
\begin{align}
\displaystyle
 |\varphi_0 - {I}_{T_0} \varphi_0 |_{W^{m,q}(T_0)}
&\leq C_*^I |T_0|^{\frac{1}{q} - \frac{1}{p}} \left( \frac{\alpha_{\max}}{\alpha_{\min}} \right)^m \left( \frac{H_{T_0}}{h_{T_0}} \right)^m h_{T_0}^{\ell+1-m} | \varphi_0 |_{W^{\ell+1,p}(T_0)}, \label{cla41}
\end{align}
where $C_*^I$ is a positive constant that is independent of $h_T$ and $H_T$, and the parameters $\alpha_{\max}$ and $\alpha_{\min}$ are defined as
\begin{align}
\displaystyle
\alpha_{\max} := \max \{ \alpha_1 , \ldots, \alpha_d \}, \quad \alpha_{\min} := \min \{ \alpha_1 , \ldots, \alpha_d \}. \label{cla42}
\end{align}

\end{th21*}

\begin{pf*}
Let $\hat{\varphi} \in W^{\ell+1,p}(\widehat{T})$. Because $0 \leq \ell \leq k$, $\mathcal{P}^{\ell} \subset \mathcal{P}^k \subset \widehat{{P}}$. Therefore, for any $\hat{\eta} \in \mathcal{P}^{\ell}$, we have $ I_{\widehat{T}} \hat{\eta} = \hat{\eta}$. Using Eqs.~\eqref{int513} and \eqref{cla40}, we obtain
\begin{align*}
\displaystyle
| \tilde{\varphi} - I_{\widehat{T}} \hat{\varphi} |_{W^{m,q}(\widehat{T})}
&\leq | \tilde{\varphi} -\hat{\eta} |_{W^{m,q}(\widehat{T})} + |  I_{\widehat{T}} ( \hat{\eta} -  \hat{\varphi} ) |_{W^{m,q}(\widehat{T})} \\
&\leq c  \| \tilde{\varphi} -\hat{\eta} \|_{W^{\ell+1,p}(\widehat{T})},
\end{align*}
where we have used the stability of the interpolation operator $I_{\widehat{T}}$; i.e.,
\begin{align*}
\displaystyle
|  I_{\widehat{T}} ( \hat{\eta} -  \hat{\varphi} ) |_{W^{m,q}(\widehat{T})} 
&\leq \sum_{i=1}^{n_0} |\hat{\chi}_i (\hat{\eta} -  \hat{\varphi})| |\hat{\theta}_i|_{W^{m,q}(\widehat{T})}  
\leq c \| \hat{\eta} - \tilde{\varphi} \|_{W^{\ell+1,p}(\widehat{T})}.
\end{align*}
Using the classic Bramble--Hilbert-type lemma (e.g., \cite[Lemma 4.3.8]{BreSco08}), we obtain
\begin{align}
\displaystyle
| \tilde{\varphi} - I_{\widehat{T}} \hat{\varphi} |_{W^{m,q}(\widehat{T})}
&\leq c \inf_{\hat{\eta} \in \mathcal{P}^{\ell}} \| \hat{\eta} - \tilde{\varphi} \|_{W^{\ell+1,p}(\widehat{T})} \leq c |\hat{\varphi}|_{W^{\ell+1,p}(\widehat{T})}. \label{cla43}
\end{align}

Inequalities \eqref{sa22}, \eqref{sa25}, \eqref{jensen}, and \eqref{cla43} yield
\begin{align}
\displaystyle
&| \varphi_0 - {I}_{T_0} \varphi_0 |_{W^{m,q}({T}_0)} \leq c | \varphi - {I}_T \varphi |_{W^{m,q}({T})} \notag \\
&\ \leq c |\det({A}_T)|^{\frac{1}{q}} \| \widetilde{A}^{-1} \|_2^m \left( \sum_{|\beta| = m}  ( \alpha^{- \beta} )^q \|  \partial^{\beta}  (\hat{\varphi} - I_{\widehat{T}} \hat{\varphi} )  \|^q_{L^q(\widehat{T})} \right)^{1/q} \notag \\
&\ \leq c |\det({A}_T)|^{\frac{1}{q}}  \| \widetilde{A}^{-1} \|_2^m \max \{ \alpha_1^{-1} , \ldots, \alpha_d^{-1} \}^{|\beta|} | \tilde{\varphi} - I_{\widehat{T}} \hat{\varphi} |_{W^{m,q}(\widehat{T})} \notag \\
&\ \leq c |\det({A}_T)|^{\frac{1}{q}} \| \widetilde{A}^{-1} \|_2^m  \alpha_{\min}^{- |\beta|}  |\hat{\varphi}|_{W^{\ell+1,p}(\widehat{T})}. \label{cla44}
\end{align}
Using inequalities \eqref{jensen} and \eqref{sa32} together with Eq.~\eqref{CN331c}, we have
\begin{align}
\displaystyle
&|\hat{\varphi}|_{W^{\ell+1,p}(\widehat{T})} \notag \\
&\ \leq  \sum_{|\gamma| = \ell + 1 -m} \sum_{|\beta| = m} \| \partial^{\beta} \partial^{\gamma} \hat{\varphi} \|_{L^{p}(\widehat{T})} \notag \\
&\ \leq c |\det({A}_T)|^{-\frac{1}{p}} \| \widetilde{{A}}  \|_2^{m}  \sum_{|\gamma| = \ell + 1 -m} \sum_{|\beta| = m} \alpha^{\beta} \sum_{|\epsilon| =  |\gamma|} {\alpha}^{\varepsilon} | \partial_{r^{(0)}}^{\epsilon} \varphi_0 |_{W^{m,p}(T_0)}  \notag \\
&\ \leq c  |\det({A}_T)|^{-\frac{1}{p}} \| \widetilde{{A}}  \|_2^{m}  \max \{ \alpha_1 , \ldots, \alpha_d \}^{|\beta|} h_{T_0}^{\ell+1-m} |\varphi_0|_{W^{\ell+1,p} (T_0)} \notag \\
&\ \leq c  |\det({A}_T)|^{-\frac{1}{p}} \| \widetilde{{A}}  \|_2^{m} \alpha_{\max}^{|\beta|}  h_{T_0}^{\ell+1-m} |\varphi_0|_{W^{\ell+1,p} (T_0)}.  \label{cla45}
\end{align}
From Eqs.~\eqref{cla44} and  \eqref{cla45} together with Eq.~\eqref{CN332}, we have the desired estimate \eqref{cla41}.
\qed
\end{pf*}

\section{Anisotropic interpolation error estimates} \label{sec_main}
%The following theorem is the one to replace \cite[Theorem 2]{IshKobTsu}.

\subsection{Main theorem}
%The following is the main theorem of this paper. It is possible to avoid using the maximum-angle condition.

Theorem A can be applied to standard isotropic elements as well as some classes of anisotropic elements. If we are concerned with anisotropic elements, it is desirable to remove the quantity $\alpha_{\max} / \alpha_{\min}$  from estimate \eqref{cla41}. To this end, we employ the approach described in \cite{Ape99} and consider the case of a finite element with $V(\widehat{T}) := \mathcal{C}({\widehat{T}})$ and $\widehat{{P}} := \mathcal{P}^{k}(\widehat{T})$ (Theorem B). However, one needs stronger assumptions to obtain the optimal estimate. When using finite elements that do not satisfy the assumptions of Theorem B (e.g., $\mathcal{P}^1$-bubble finite element), we have to use Theorem A. In these cases, it may not be possible to obtain optimal order estimates if the shape-regularity condition is violated.

\begin{thrBA*} \label{thA}
Let $\{ \widehat{T} , \widehat{{P}} , \widehat{\Sigma} \}$ be a finite element with the normed vector space $V(\widehat{T}) := \mathcal{C}({\widehat{T}})$ and $\widehat{{P}} := \mathcal{P}^{k}(\widehat{T})$ with $k \geq 1$.  Let $I_{\widehat{T}} : V({\widehat{T}}) \to   \widehat{{P}}$ be a linear operator. Fix $\ell \in \mathbb{N}$, $m \in \mathbb{N}_0$, and $p,q \in [1,\infty]$ such that $0 \leq m \leq \ell  \leq k+1$, $\ell - m \geq 1$, and assume that the embeddings \eqref{emmed} and \eqref{emmed2} with $s := \ell - m$ hold. Let $\beta$ be a multi-index with $|\beta| = m$. We set $j := \dim (\partial^{\beta} \mathcal{P}^{k})$. Assume that there exist linear functionals $\mathscr{F}_i$, $i=1,\ldots,j$, such that
\begin{align}
\displaystyle
&\mathscr{F}_i \in W^{\ell  - m ,p} (\widehat{T})^{\prime},\quad \forall i = 1,\ldots,j, \label{main39}\\
&\mathscr{F}_i ( \partial^{\beta} (\hat{\varphi} - I_{\widehat{T}} \hat{\varphi}) ) = 0 \quad  \forall i = 1,\ldots,j, \quad \forall \hat{\varphi} \in \mathcal{C}({\widehat{T}}): \ \partial^{\beta} \hat{\varphi} \in W^{\ell - m , p} (\widehat{T}), \label{main40} \\
&\hat{\eta} \in \mathcal{P}^{k}, \quad \mathscr{F}_i(\partial^{\beta} \hat{\eta}) = 0 \quad \forall i = 1,\ldots,j \quad \Rightarrow \quad \partial^{\beta} \hat{\eta} = 0. \label{main41}
\end{align}
For any $\hat{\varphi} \in W^{\ell   , p} (\widehat{T}) \cap \mathcal{C}({\widehat{T}})$, we set ${\varphi}_0 := \hat{\varphi} \circ {\Phi}^{-1}$. If Assumption \ref{ass1} is imposed, it holds that
\begin{align}
\displaystyle
&| {\varphi}_0 - I_{{T}_0} {\varphi}_0 |_{W^{m,q}({T}_0)} \nonumber \\
&\quad \leq  C_1^{TB} |T_0|^{\frac{1}{q} - \frac{1}{p}} \left( \frac{H_{T_0}}{h_{T_0}} \right)^m \sum_{|\gamma| = \ell-m} \mathscr{H}^{\gamma}  | \partial^{\gamma}  {(\varphi_0 \circ \Phi_{T_0})} |_{W^{ m ,p}(\Phi_{T_0}^{-1}(T_0))}, \label{main42}
\end{align}
where $C_1^{TB}$ is a positive constant that is independent of $h_{T_0}$ and $H_{T_0}$. Furthermore,  if Assumption \ref{ass1} is not imposed, it holds that
\begin{align}
\displaystyle
&| {\varphi}_0 - I_{{T}_0} {\varphi}_0 |_{W^{m,q}({T}_0)} \nonumber \\
&\quad \leq  C_2^{TB} |T_0|^{\frac{1}{q} - \frac{1}{p}} \left( \frac{H_{T_0}}{h_{T_0}} \right)^m \sum_{|\gamma| = \ell-m} {\alpha}^{\gamma}  | \partial_{r^{(0)}}^{\gamma}  {\varphi_0} |_{W^{ m ,p}(T_0)}, \label{main43}
\end{align}
where $C_2^{TB}$ is a positive constant that is independent of $h_{T_0}$ and $H_{T_0}$. 
\end{thrBA*}

\begin{pf*}
The introduction of the functionals $\mathscr{F}_i$ follows from \cite{Ape99}. In fact, under the same assumptions as made in Theorem B, we have (see \cite[Lemma 2.2]{Ape99})
\begin{align}
\displaystyle
\| \partial^{\beta} (\hat{\varphi} - I_{\widehat{T}} \hat{\varphi}) \|_{L^q(\widehat{T})}
&\leq C^B | \partial^{\beta}  \hat{\varphi} |_{W^{\ell  - m ,p}(\widehat{T})}, \label{main44}
\end{align}
where $|\beta| = m$, $\hat{\varphi} \in \mathcal{C}({\widehat{T}})$, and $\partial^{\beta} \hat{\varphi} \in W^{\ell - m , p} (\widehat{T})$.

Inequalities \eqref{sa22}, \eqref{sa25}, \eqref{jensen}, and \eqref{main44} yield
\begin{align}
\displaystyle
&| \varphi_0 - {I}_{T_0} \varphi_0 |_{W^{m,q}({T}_0)} \leq c | \varphi - {I}_T \varphi |_{W^{m,q}({T})} \notag \\
&\ \leq c |\det({A}_T)|^{\frac{1}{q}} \| \widetilde{A}^{-1} \|_2^m \left( \sum_{|\beta| = m}  ( \alpha^{- \beta} )^q \|  \partial^{\beta}  (\hat{\varphi} - I_{\widehat{T}} \hat{\varphi} )  \|^q_{L^q(\widehat{T})} \right)^{1/q} \notag \\
&\ \leq c |\det({A}_T)|^{\frac{1}{q}}  \| \widetilde{A}^{-1} \|_2^m  \sum_{|\beta| = m}  ( \alpha^{- \beta} ) \|  \partial^{\beta}  (\hat{\varphi} - I_{\widehat{T}} \hat{\varphi} )  \|_{L^q(\widehat{T})} \notag \\
&\ \leq c |\det({A}_T)|^{\frac{1}{q}} \| \widetilde{A}^{-1} \|_2^m  \sum_{|\beta| = m}  ( \alpha^{- \beta} )  | \partial^{\beta}  \hat{\varphi} |_{W^{\ell  - m ,p}(\widehat{T})}. \label{main45}
\end{align}

If Assumption \ref{ass1} is imposed, then using inequalities \eqref{jensen} and \eqref{sa29} leads to
\begin{align}
\displaystyle
& \sum_{|\beta| = m}  ( \alpha^{- \beta} )  | \partial^{\beta}  \hat{\varphi} |_{W^{\ell  - m ,p}(\widehat{T})} \notag \\
 &\quad \leq  \sum_{|\gamma| = \ell-m} \sum_{|\beta| = m}  ( \alpha^{- \beta} )  \| \partial^{\beta + \gamma}  \hat{\varphi} \|_{L^{p}(\widehat{T})} \notag \\
 &\quad \leq  c |\det({A}_T)|^{-\frac{1}{p}} \| \widetilde{A} \|^m_2 \sum_{|\gamma| = \ell-m} \sum_{|\beta| = m}  ( \alpha^{- \beta} )  \alpha^{\beta} \sum_{|\epsilon| =  |\gamma|} \mathscr{H}^{\varepsilon} | \partial^{\epsilon} \varphi |_{W^{m,p}(T)} \notag \\
 &\quad \leq c  |\det({A}_T)|^{-\frac{1}{p}} \| \widetilde{A} \|^m_2 \sum_{|\epsilon| =  \ell-m} \mathscr{H}^{\varepsilon} | \partial^{\epsilon} \varphi |_{W^{m,p}(T)}. \label{main46}
\end{align}
From Eqs.~\eqref{main45} and  \eqref{main46} together with Eqs.~\eqref{CN331} and \eqref{CN332}, we have the desired estimate \eqref{main42} using $T = \Phi_{T_0}^{-1}(T_0)$ and $\varphi = \varphi_0 \circ \Phi_{T_0}$.

If Assumption \ref{ass1} is not imposed, then an analogous argument using inequality \eqref{sa32} instead of \eqref{sa29} yields estimate \eqref{main43}.
\qed
\end{pf*}

\begin{Ex}
Specific finite elements satisfying conditions \eqref{main39}, \eqref{main40}, and \eqref{main41} are given in \cite{ApeDob92} and \cite{Ape99}; see also Section \ref{sec71}.
\end{Ex}

\begin{Rem} \label{counter_ex}
Finite elements that do not satisfy conditions \eqref{main39}, \eqref{main40}, and \eqref{main41} can be found in \cite[Table 3]{ApeDob92}; e.g., the $\mathcal{P}^1$-bubble finite element and the $\mathcal{P}^3$ Hermite finite element. In these cases, Theorem A can be applied.
\end{Rem}

\subsection{Examples satisfying conditions \eqref{main39}, \eqref{main40}, and \eqref{main41} in Theorem B} \label{sec71}

\begin{coro} \label{coro1}
Let $\{ \widehat{T} , \widehat{{P}} , \widehat{\Sigma} \}$ be the Lagrange finite element with  $V(\widehat{T}) := \mathcal{C}({\widehat{T}})$ and $\widehat{{P}} := \mathcal{P}^{k}(\widehat{T})$ for $k \geq 1$. Let $I_T: V({T})  \to {{P}}$ be the corresponding local Lagrange interpolation operator. Let  $m \in \mathbb{N}_0$, $\ell \in \mathbb{N}$, and $p \in \mathbb{R}$ be such that $0 \leq m \leq \ell \leq k+1$ and
\begin{align*}
\displaystyle
&d=2: \ 
 \begin{cases}
p \in (2,\infty] \quad \text{if $m=0$, $\ell = 1$},\\
p \in [1,\infty] \quad \text{if $m=0$, $\ell \geq 2$ or $m \geq 1$, $\ell - m \geq 1$},
\end{cases} \\
&d=3: \ \begin{cases}
p \in \left(\frac{3}{\ell}, \infty \right] \quad \text{if $m=0$, $\ell=1,2$},\\
p \in (2,\infty] \quad \text{if $m \geq 1$,  $\ell-m = 1$}, \\
p \in [1,\infty] \quad \text{if $m=0$, $\ell \geq 3$ or $m \geq 1$,  $\ell-m \geq 2$}.
\end{cases}
\end{align*}
We set $q \in [1,\infty]$ such that $W^{\ell-m,p}(\widehat{T}) \hookrightarrow L^q(\widehat{T})$. Then, for all $\hat{\varphi} \in W^{\ell ,p}(\widehat{T})$ with ${\varphi}_0 := \hat{\varphi} \circ {\Phi}^{-1}$, we recover Eq.~\eqref{main42} if Assumption \ref{ass1} is imposed, and Eq.~\eqref{main43} holds.

Furthermore, for any $\hat{\varphi} \in \mathcal{C}(\widehat{T})$ with ${\varphi}_0 := \hat{\varphi} \circ {\Phi}^{-1}$, it holds that
\begin{align*}
\displaystyle
\| {\varphi}_0 - I_{{T}_0} {\varphi}_0 \|_{L^{\infty}({T}_0)}
\leq  c  \|  {\varphi}_0 \|_{L^{\infty}(T_0)}.
\end{align*}
\end{coro}

\begin{pf*}
The existence of functionals satisfying Eqs.~\eqref{main39}, \eqref{main40}, and \eqref{main41} is shown in the proof of \cite[Lemma 2.4]{Ape99} for $d=2$ and in the proof of \cite[Lemma 2.6]{Ape99} for $d=3$. Inequality \eqref{main44} then holds. This implies that estimates \eqref{main42} and \eqref{main43} hold.
\qed
\end{pf*}

Setting $V(T) := \mathcal{C}({{T}})$, we define the nodal Crouzeix--Raviart interpolation operators as
\begin{align*}
\displaystyle
I_{T}^{CR,S}: V(T) \ni \varphi  \mapsto I_{T}^{CR,S} \varphi := \sum_{i=1}^{d+1} \varphi (x_{F_i}) \theta_i \in \mathcal{P}^1.
\end{align*}

\begin{coro} \label{coro2}
Let $\{ \widehat{T} , \widehat{{P}} , \widehat{\Sigma} \}$ be the Crouzeix--Raviart finite element with  $V(\widehat{T}) := \mathcal{C}({\widehat{T}})$ and $\widehat{{P}} := \mathcal{P}^{1}(\widehat{T})$. Set $I_T ;= I_{T}^{CR,S}$. Let  $m \in \mathbb{N}_0$, $\ell \in \mathbb{N}$, and $p \in \mathbb{R}$ be such that
\begin{align*}
\displaystyle
&d=2: \ 
 \begin{cases}
p \in (2,\infty]  \quad \text{if $m=0$, $\ell = 1$},\\
p \in [1,\infty] \quad \text{if $m=0$, $\ell = 2$ or $m=1$, $\ell = 2$},
\end{cases} \\
&d=3: \ \begin{cases}
%\displaystyle
p \in \left(\frac{3}{\ell} , \infty \right] \quad \text{if $m=0$, $\ell=1,2$},\\
\displaystyle
p \in (2,\infty] \quad \text{if $m = 1$,  $\ell = 2$}.
\end{cases}
\end{align*}
Set $q \in [1,\infty]$ such that $W^{\ell-m,p}(\widehat{T}) \hookrightarrow L^q(\widehat{T})$. Then, for all $\hat{\varphi} \in W^{\ell ,p}(\widehat{T})$ with ${\varphi}_0 := \hat{\varphi} \circ {\Phi}^{-1}$, we recover Eq.~\eqref{main42} if Assumption \ref{ass1} is imposed, and Eq.~\eqref{main43} holds.

Furthermore, for any $\hat{\varphi} \in \mathcal{C}(\widehat{T})$ with ${\varphi}_0 := \hat{\varphi} \circ {\Phi}^{-1}$, it holds that
\begin{align*}
\displaystyle
\| {\varphi}_0 - I_{{T}_0} {\varphi}_0 \|_{L^{\infty}({T}_0)}
\leq  c  \|  {\varphi}_0 \|_{L^{\infty}(T_0)}.
\end{align*}
\end{coro}

\begin{pf*}
%For $k=1$, we show that there exists functionals $\mathscr{F}_i$ satisfying \eqref{4}, \eqref{5}, and \eqref{6} in Theorem B for each $\ell$ and $m$.
For $k=1$, we only introduce functionals $\mathscr{F}_i$ satisfying Eqs.~\eqref{main39}, \eqref{main40}, and \eqref{main41} in Theorem B for each $\ell$ and $m$.

%Let $m=0$. From the Sobolev embedding theorem (e.g., see \cite[Theorem 4.12]{AdaFou03}, \cite[pp. 27, 28]{Ape99}, \cite[p. 114]{Cia02}, \cite[Corollary B.43]{ErnGue04}, \cite[Theorem 1.3]{GirRav86}), we have $W^{\ell,p}(\widehat{T}) \subset \mathcal{C}^0(\widehat{T})$ with $1 \leq p \leq \infty$ and $d \< \ell p$. Under this condition, we use
Let $m=0$. From the Sobolev embedding theorem, we have $W^{\ell,p}(\widehat{T}) \subset \mathcal{C}^0(\widehat{T})$ with $1 \< p \leq \infty$, $d \< \ell p$ or $p=1$, $d \leq \ell$. Under this condition, we use
\begin{align*}
\displaystyle
\mathscr{F}_i(\hat{\varphi}) := \hat{\varphi} (\hat{x}_{\widehat{F}_i}), \quad \hat{\varphi} \in W^{\ell,p}(\widehat{T}), \quad i=1,\ldots, d+1.
\end{align*}

Let $d=2$ and $m=1$ ($\ell = 2$). We set $\beta = (1,0)$. Then, we have that $j = \dim( \partial^{\beta} \mathcal{P}^1 ) = 1$. We consider a functional
\begin{align*}
\displaystyle
\mathscr{F}_1(\hat{\varphi}) := \int_{0}^{\frac{1}{2}} \hat{\varphi}(\hat{x}_1,1/2) d \hat{x}_1, \quad \hat{\varphi} \in W^{2,p}(\widehat{T}).
\end{align*}
By an analogous argument, we can set a functional for the case $\beta = (0,1)$.

Let $d=3$ and $m=1$ ($\ell = 2$). We first consider Type (\roman{sone}) in Section \ref{reference3d}. That is, the reference element is $\widehat{T} = \conv \{ 0,e_1, e_2,e_3 \}$. Here, $e_1, \ldots, e_3 \in \mathbb{R}^3$ form the canonical basis. We set $\beta = (1,0,0)$ and consider the functional
\begin{align*}
\displaystyle
\mathscr{F}_1(\hat{\varphi}) := \int_{0}^{\frac{1}{3}} \hat{\varphi}(\hat{x}_1,1/3,1/3) d \hat{x}_1, \quad \hat{\varphi} \in W^{2,p}(\widehat{T}).
\end{align*}
We now consider Type (\roman{stwo}) in Section \ref{reference3d}. That is, the reference element is $\widehat{T} = \conv \{ 0,e_1, e_1 + e_2 , e_3 \}$. We set $\beta = (1,0,0)$ and consider the functional
\begin{align*}
\displaystyle
\Phi_1(\hat{\varphi}) := \int_{\frac{1}{3}}^{\frac{2}{3}} \hat{\varphi}(\hat{x}_1,1/3,1/3) d \hat{x}_1, \quad \hat{\varphi} \in W^{2,p}(\widehat{T}).
\end{align*}
By an analogous argument, we can set functionals for cases $\beta = (0,1,0), (0,0,1)$.

When $m = \ell = 0$ and $p = \infty$, we can easily check that
\begin{align*}
\displaystyle
\| \hat{\varphi} - I_{\widehat{T}}^{CR,S} \hat{\varphi} \|_{L^{\infty}(\widehat{T})}
&\leq c \| \hat{\varphi} \|_{L^{\infty}(\widehat{T})}.
\end{align*}
\begin{comment}
because we have
\begin{align*}
\displaystyle
| ( I_{\widehat{T}}^{CR,S}  \hat{\varphi} ) (\hat{x}) |
&\leq \sum_{i=1}^{d+1} | \hat{\varphi} (\hat{x}_{\widehat{F}_i})| | \hat{\theta}_i (\hat{x})|
\leq (d+1) \left( \max_{1 \leq i \leq d+1} \| \hat{\theta}_i \|_{L^{\infty}(\widehat{T})} \right) \| \hat{\varphi} \|_{L^{\infty}(\widehat{T})}.
\end{align*}
\end{comment}
\qed
\end{pf*}

\section{Concluding remarks}
As our concluding remarks, we identify several topics related to the results described in this paper. 

\subsection{Good elements or not for $d=2,3$?} \label{goodbad}
In this subsection, we consider good elements on meshes. Here, we define ``good elements" on meshes as those for which there exists some $\gamma_0 \> 0$ satisfying Eq.~\eqref{NewGeo}. We treat a ``Right-angled triangle," ``Blade," and ``Dagger" for $d=2$, and a ``Spire," ``Spear," ``Spindle," ``Spike," ``Splinter," and "Sliver" for $d=3$, as introduced in \cite{Cheetal00}. We present the quantities $\alpha_{\max} / \alpha_{\min}$ and $H_{T_0}/h_{T_0}$ for these elements.

\subsubsection{Isotropic mesh} \label{Isomesh}
We consider the following condition. There exists a constant $\gamma_1 \> 0$ such that, for any $\mathbb{T}_h \in \{ \mathbb{T}_h \}$ and any simplex $T_0 \in \mathbb{T}_h$, we have
\begin{align}
\displaystyle
|T_0| \geq \gamma_1 h_{T_0}^d. \label{geo3}
\end{align}
Condition \eqref{geo3} is equivalent to the shape-regularity condition; see \cite[Theorem 1]{BraKorKri08}.

%\begin{Ex} \label{iso}
If geometric condition \eqref{geo3} is satisfied, it holds that
\begin{align*}
\displaystyle
\frac{H_{T_0}}{h_{T_0}} \leq \frac{h_{T_0}^d}{|T_0|} \leq \frac{1}{\gamma_1}, \quad \frac{\alpha_{\max}}{\alpha_{\min}} \leq  c \frac{h_{T_0}}{\alpha_2} \leq c \frac{h_{T_0}^d}{|T_0|} \leq \frac{c}{\gamma_1}.
%\frac{\alpha_{\max}}{\alpha_{\min}} \leq c \frac{h_T^d}{|T|} = c \frac{h_{T_0}^d}{|T_0|} \leq \frac{c}{\gamma_1}.
\end{align*}
If  $p=q$ in Theorem A, one can obtain the optimal order $h_{T_0}^{\ell+1-m}$. In this case, elements satisfying geometric condition \eqref{geo3} are ``good."
%\end{Ex}

\subsubsection{Anisotropic mesh: two-dimensional case}
Let $S_0 \subset \mathbb{R}^2$ be a triangle. Let $0 \< s \ll 1$, $s \in \mathbb{R}$, and $\varepsilon,\delta,\gamma \in \mathbb{R}$. A dagger has one short edge and a blade has no short edge.

%All examples degenerate in the $y$-axis direction. 

\begin{Ex}[Right-angled triangle] \label{Note2}
Let $S_0 \subset \mathbb{R}^2$ be the simplex with vertices $x_1 := (0,0)^T$, $x_2 := (s,0)^T$, and $x_3 := (0,s^{\varepsilon})^T$ with $1 \< \varepsilon$. Then, we have that $\alpha_1 = s$ and $\alpha_2 =  s^{\varepsilon}$; i.e.,
\begin{align*}
\displaystyle
\frac{\alpha_{\max}}{\alpha_{\min}} &\leq s^{1 - \varepsilon} \to \infty \quad \text{as $s \to 0$}, \quad \frac{H_{S_0}}{h_{S_0}} =2.
\end{align*}
In this case, the element $S_0$ is ``good."
\end{Ex}

\begin{Ex}[Dagger] \label{Note1}
Let $S_0 \subset \mathbb{R}^2$ be the simplex with vertices $x_1 := (0,0)^T$, $x_2 := (s,0)^T$, and $x_3 := (s^{\delta},s^{\varepsilon})^T$ with $1 \< \varepsilon \< \delta $. Then, we have that $\alpha_1 = \sqrt{(s - s^{\delta})^2 + s^{2 \varepsilon}}$ and $\alpha_2 = \sqrt{s^{2 \delta} + s^{2 \varepsilon}}$; i.e.,
\begin{align*}
\displaystyle
\frac{\alpha_{\max}}{\alpha_{\min}} &= \frac{ \sqrt{(s - s^{\delta})^2 + s^{2 \varepsilon}}}{\sqrt{s^{2 \delta} + s^{2 \varepsilon}}} \leq c s^{1 - \varepsilon} \to \infty \quad \text{as $s \to 0$},\\
\frac{H_{S_0}}{h_{S_0}} &= \frac{ \sqrt{(s- s^{\delta})^2 + s^{2 \varepsilon}} \sqrt{s^{2 \delta} + s^{2 \varepsilon}}}{\frac{1}{2} s^{1 + \varepsilon}} \leq c.
\end{align*}
In this case, the element $S_0$ is ``good."
\end{Ex}

\begin{Rem}
In the above examples, $\alpha_2 \approx \alpha_2 t = \mathscr{H}_2$ holds. That is, the good element $S_0 \subset \mathbb{R}^2$ satisfies conditions such as $\alpha_2 \approx \alpha_2 t = \mathscr{H}_2$.
\end{Rem}

\begin{Ex}[Blade] \label{Ex4}
Let $S_0 \subset \mathbb{R}^2$ be the simplex with vertices $x_1 := (0,0)^T$, $x_2 := (2s,0)^T$, and $x_3 := (s ,s^{\varepsilon})^T$ with $1 \< \varepsilon $. Then, we have that $\alpha_1 = \alpha_2 = \sqrt{s^{2} + s^{2 \varepsilon}}$; i.e.,
\begin{align*}
\displaystyle
\frac{\alpha_{\max}}{\alpha_{\min}} = 1, \quad \frac{H_{S_0}}{h_{S_0}} = \frac{s^{2} + s^{2 \varepsilon}}{s^{1 + \varepsilon}} \to \infty \quad \text{as $s \to 0$}.
\end{align*}
In this case, the element $S_0$ is ``not good."
\end{Ex}

\begin{Ex}[Dagger] \label{Ex5}
Let $S_0 \subset \mathbb{R}^2$ be the simplex with vertices $x_1 := (0,0)^T$, $x_2 := (s,0)^T$, and $x_3 := (s^{\delta},s^{\varepsilon})^T$ with $1 \< \delta \< \varepsilon $. Then, we have that $\alpha_1 = \sqrt{(s - s^{\delta})^2 + s^{2 \varepsilon}}$ and $\alpha_2 = \sqrt{s^{2 \delta} + s^{2 \varepsilon}}$; i.e.,
\begin{align*}
\displaystyle
\frac{\alpha_{\max}}{\alpha_{\min}} &= \frac{ \sqrt{(s - s^{\delta})^2 + s^{2 \varepsilon}}}{\sqrt{s^{2 \delta} + s^{2 \varepsilon}}} \leq c s^{1 - \delta} \to \infty \quad \text{as $s \to 0$},\\
\frac{H_{S_0}}{h_{S_0}} &= \frac{ \sqrt{(s- s^{\delta})^2 + s^{2 \varepsilon}} \sqrt{s^{2 \delta} + s^{2 \varepsilon}}}{\frac{1}{2} s^{1 + \varepsilon}} \leq c s^{\delta - \varepsilon} \to \infty \quad \text{as $s \to 0$}.
\end{align*}
In this case, the element $S_0$ is ``not good."
\end{Ex}

%\begin{Rem}
%The above two examples indicate that if the condition \eqref{NewGeo} is violated, the convergence order in both Theorem A and Theorem B deteriorates even if $p=q$. In some cases, the interpolation error estimates diverge as meshes get finer.
%\end{Rem}

\begin{comment}
\begin{figure}[tbhp]
\vspace{-5cm}
  \includegraphics[bb=0 0 658 448,scale=0.55]{rat.png}
\caption{Example 2: Right-angled triangle}
\label{rat}
\end{figure}

\begin{figure}[tbhp]
\vspace{-5cm}
  \includegraphics[bb=0 0 658 448,scale=0.55]{dagger.png}
\caption{Examples 3 and 5: Dagger}
\label{rat}
\end{figure}

\begin{figure}[htbp]
\vspace{-4cm}
  \includegraphics[bb=0 0 658 448,scale=0.55]{blade.png}
\caption{Example 4: Blade}
\label{rat}
\end{figure}
\end{comment}

\begin{figure}[htbp]
\vspace{-4cm}
  \begin{minipage}[b]{0.52\linewidth}
    \centering
   \includegraphics[bb=0 0 658 448,scale=0.55]{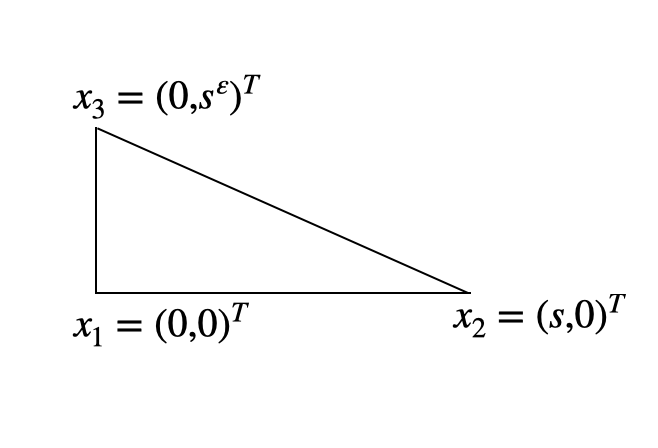}
    \caption{Example 2: Right-angled triangle}
    \label{RAT}
  \end{minipage}
  \begin{minipage}[b]{0.52\linewidth}
    \centering
     \includegraphics[bb=0 0 658 448,scale=0.55]{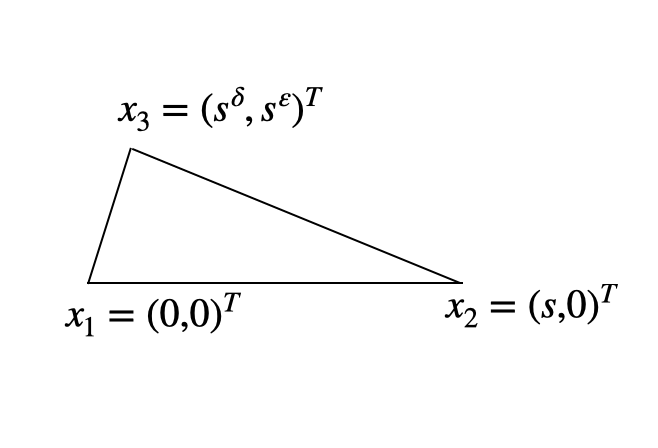}
    \caption{Examples 3 and 5: Dagger}
  \end{minipage}
\end{figure}

\begin{figure}[htbp]
\vspace{-4cm}
  \includegraphics[bb=0 0 658 448,scale=0.55]{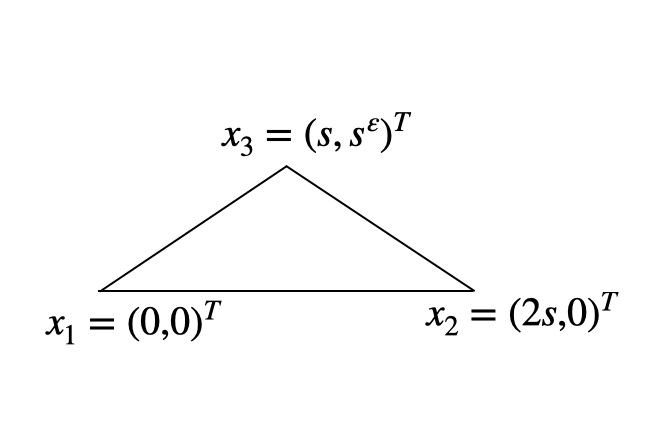}
\caption{Example 4: Blade}
\label{rat}
\end{figure}

\subsubsection{Anisotropic mesh: three-dimensional case}
\begin{Ex}
%Let $T_0 \subset \mathbb{R}^3$ be a tetrahedron and let $T \subset \mathbb{R}^3$ the tetrahedron satisfying Condition \ref{cond2} through appropriate rotation, translation, and mirror imaging. . 
Let $T_0 \subset \mathbb{R}^3$ be a tetrahedron. Let $S_0$ be the base of $T_0$; i.e., $S_0 = \triangle x_1 x_2 x_3$. Recall that
\begin{align}
\displaystyle
\frac{H_{T_0}}{h_{T_0}} = \frac{\alpha_1 \alpha_2 \alpha_3}{|T_0|} = \frac{\alpha_1 \alpha_2}{\frac{1}{2} \alpha_1 \alpha_2 t_1} \frac{\alpha_3}{\frac{1}{3} \alpha_3 t_2} \leq  \frac{H_{S_0}}{h_{S_0}}  \frac{\alpha_3}{\frac{1}{3} \mathscr{H}_3}. \label{bad56}
\end{align}
If the triangle $S_0$ is ``not good," such as in Examples \ref{Ex4} and \ref{Ex5}, the quantity in Eq.~\eqref{bad56} may diverge. In the following, we consider the case in which the triangle $S_0$ is ``good."

Assume that there exists a positive constant $M$ such that $\frac{H_{S_0}}{h_{S_0}} \leq M$. For simplicity, we set $x_1 := (0,0,0)^T$, $x_2 := (2s,0,0)^T$, and $x_3 := (2s - \sqrt{4 s^2 - s^{2 \gamma}}, s^{\gamma},0)^T$ with $1 \< \gamma$. Then,
\begin{align*}
\displaystyle
\alpha_1 = 2s, \ \alpha_2 = \sqrt{\frac{4 s^{2 \gamma}}{2 + \sqrt{4  - s^{2 \gamma -2} }}},
\end{align*}
and because $\alpha_{\max} \approx c s$,
\begin{align*}
\displaystyle
\frac{\alpha_{\max}}{\alpha_{\min}} &\leq \frac{c s}{\alpha_2} \leq c s^{1 - \gamma} \to \infty \quad \text{as $s \to 0$}.
\end{align*}
%Then, because $\alpha_2$ is the minimum edge of $T$,
%\begin{align*}
%\displaystyle
%\frac{\alpha_{\max}}{\alpha_{\min}} \leq c s^{1 - \gamma} \to \infty \quad \text{as $s \to 0$}.
%\end{align*}

If we set $x_4 := (s,0,s^{\varepsilon})^T$ with $1 \< \varepsilon$, the triangle $\triangle x_1 x_2 x_4$ is the blade (Example \ref{Ex4}). Then,
\begin{align*}
\displaystyle
\alpha_3 = \sqrt{s^2 + s^{2 \varepsilon}}.
\end{align*}
Thus, we have
\begin{align*}
\displaystyle
\frac{H_{T_0}}{h_{T_0}} \leq c \frac{s^{2s + \gamma}}{s^{1+\gamma+\varepsilon}} \leq c s^{1-\varepsilon} \to \infty \quad \text{as $s \to 0$}.
\end{align*}
In this case, the element $T_0$ is ``not good."

If we set $x_4 := (s^{\delta},0,s^{\varepsilon})^T$ with $1 \< \delta \< \varepsilon \< \gamma$, the triangle $\triangle x_1 x_2 x_4$ is the dagger (Example \ref{Ex5}). Then,
\begin{align*}
\displaystyle
\alpha_3 = \sqrt{s^{2 \delta} + s^{2 \varepsilon}}.
\end{align*}
Thus, we have
\begin{align*}
\displaystyle
\frac{H_{T_0}}{h_{T_0}} \leq c \frac{s^{1+\gamma + \delta}}{s^{1+\gamma+\varepsilon}} \leq c s^{\delta - \varepsilon} \to \infty \quad \text{as $s \to 0$}.
\end{align*}
In this case, the element $T_0$ is ``not good."

If we set $x_4 := (s^{\delta},0,s^{\varepsilon})^T$ with $1  \< \varepsilon \< \delta \< \gamma$, the triangle $\triangle x_1 x_2 x_4$ is the dagger (Example \ref{Note1}). Then,
\begin{align*}
\displaystyle
\alpha_3 = \sqrt{s^{2 \delta} + s^{2 \varepsilon}}.
\end{align*}
Thus, we have
\begin{align*}
\displaystyle
\frac{H_{T_0}}{h_{T_0}} \leq c \frac{s^{1+\gamma + \varepsilon} }{ s^{1+\gamma+\varepsilon}} \leq c.
\end{align*}
In this case, the element $T_0$ is ``good" and  $\alpha_3 \approx \alpha_3 t_2 = \mathscr{H}_3$ holds.

\begin{figure}[htbp]
\vspace{-5cm}
  \includegraphics[bb=0 0 774 530,scale=0.55]{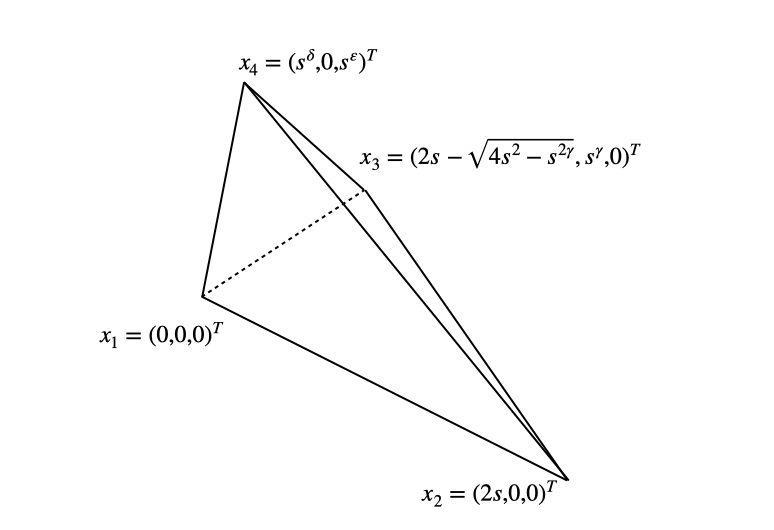}
\caption{Example 6}
\label{rat}
\end{figure}
\end{Ex}

\begin{Ex}
In \cite{Cheetal00}, the spire has a cycle of three daggers among its four triangles. The splinter has four daggers. The spear and spike have two daggers and two blades as triangles. The spindle has four blades as triangles.

\begin{figure}[htbp]
\vspace{-2cm}
  \begin{minipage}[b]{0.3\linewidth}
    \centering
    \includegraphics[bb=0 0 520 320,scale=0.35]{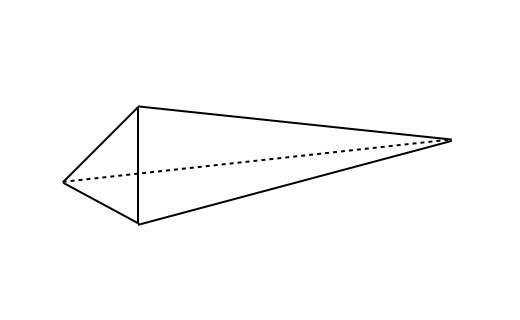}
    \caption{Spire}
  \end{minipage}
  \begin{minipage}[b]{0.3\linewidth}
    \centering
     \includegraphics[bb=0 0 520 320,scale=0.35]{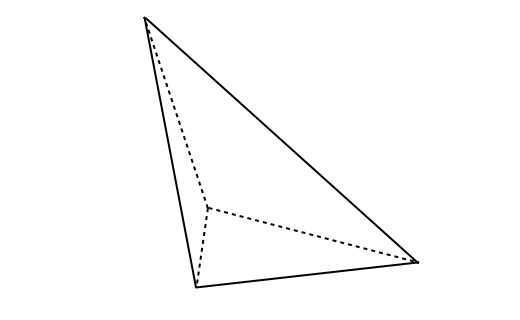}
    \caption{Spear}
  \end{minipage}
  \begin{minipage}[b]{0.3\linewidth}
    \centering
     \includegraphics[bb=0 0 520 320,scale=0.35]{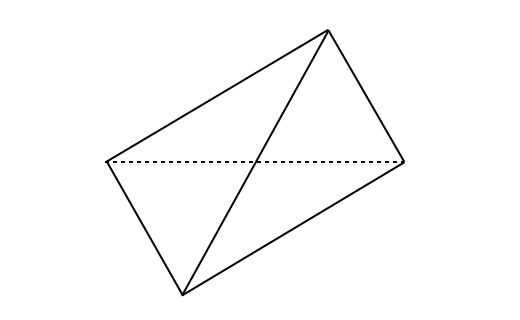}
    \caption{Spindle}
  \end{minipage}
\end{figure}

\begin{figure}[htbp]
\vspace{-2cm}
  \begin{minipage}[b]{0.3\linewidth}
    \centering
    \includegraphics[bb=0 0 520 320,scale=0.35]{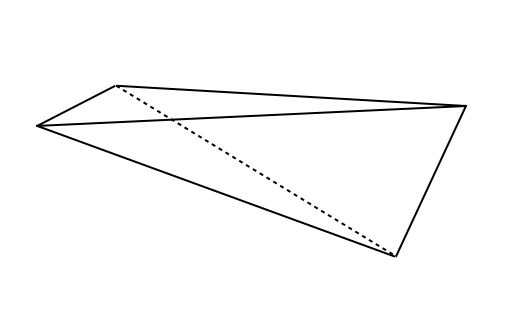}
    \caption{Spike}
  \end{minipage}
  \begin{minipage}[b]{0.3\linewidth}
    \centering
     \includegraphics[bb=0 0 520 320,scale=0.35]{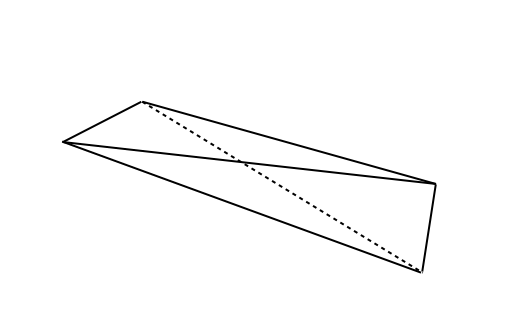}
    \caption{Splinter}
  \end{minipage}
\end{figure}

\end{Ex}

%\begin{Ex}[Spire, Spear, Spindle, Spike, Splinter]
%In \cite{CheShiZha04}, the spire has a cycle of three daggers among its four triangles. The splinter has four daggers. The spear and the spike have two daggers and two blades as triangles. The spindle has four blades as triangles.
%\end{Ex}

\begin{Rem}
The above examples reveal that the good element $T_0 \subset \mathbb{R}^3$ satisfies conditions such as $\alpha_2 \approx \alpha_2 t_1 = \mathscr{H}_2$ and $\alpha_3 \approx \alpha_3 t_2 = \mathscr{H}_3$.
\end{Rem}

\begin{Ex} \label{Ex_sliver}
Using an element $T_0$ called a \textit{Sliver}, we compare the three quantities $\frac{h_{T_0}^3}{|{T_0}|}$, $\frac{H_{T_0}}{h_{T_0}}$, and $\frac{R_3}{h_{T_0}}$, where the parameter $R_3$ denotes the circumradius of $T_0$.

Let $T_0 \subset \mathbb{R}^3$ be the simplex with vertices $x_1 := (s^{\varepsilon_2},0,0)^T$, $x_2 := (-s^{\varepsilon_2},0,0)^T$, $x_3 := (0,-s,s^{\varepsilon_1})^T$, and $x_4 := (0,s,s^{\varepsilon_1})^T$ ($\varepsilon_1, \varepsilon_2 \> 1$), where $s := \frac{1}{N}$, $N \in \mathbb{N}$. Let $L_i$ ($1 \leq i \leq 6$) be the edges of $T_0$ with $\alpha_{\min} = L_1 \leq L_2 \leq \cdots \leq L_6 = h_{T_0}$. Recall that $\alpha_{\max} \approx h_{T_0}$ and
\begin{align*}
\displaystyle
\frac{\alpha_{\max}}{\alpha_{\min}} \leq c \frac{L_6}{L_1}, \quad \frac{H_{T_0}}{h_{T_0}} = \frac{L_1 L_2}{|T_0|} h_{T_0}.
\end{align*}

\begin{table}[htbp]
\caption{$h_{T_0}^3/{|T_0|}$, $H_{T_0}/h_{T_0}$, and $R_3/h_{T_0}$ ($\varepsilon_1 = 1.5$, $\varepsilon_2 = 1.0$)}
\centering
\begin{tabular}{l | l | l | l | l | l } \hline
$N$ &  $s$ & $L_6 / L_1$ & $h_{T_0}^3/{|T_0|}$  & $H_{T_0}/h_{T_0}$ & $R_3/h_{T_0}$\\ \hline \hline
%64 & 1.5625e-02 & 96  &  9.6750e+01  & 1.0020e+00    \\
32 & 3.1250e-02 & 1.4033 & 6.7882e+01   & 3.4471e+01  &  5.0195e-01     \\
64 & 1.5625e-02  & 1.4087 & 9.6000e+01  & 4.8375e+01  &  5.0098e-01  \\
128 &7.8125e-03 &  1.4115 & 1.3576e+02  &6.8147e+01 &  5.0049e-01   \\
%1024 & 9.7656e-04  &   &     \\
\hline
\end{tabular}
\label{sliver1}
\end{table}

\begin{table}[ht]
\caption{$h_{T_0}^3/{|T_0|}$, $H_{T_0}/h_{T_0}$, and $R_3/h_{T_0}$ ($\varepsilon_1 = 1.0$, $\varepsilon_2 = 1.5$)}
\centering
\begin{tabular}{l | l | l | l | l | l } \hline
$N$ &  $s$ & $L_6 / L_1$ & $h_{T_0}^3/{|T_0|}$  & $H_{T_0}/h_{T_0}$ & $R_3/h_{T_0}$\\ \hline \hline
%64 & 1.5625e-02 & 96 & 1.2094e+01   &  1.0000e+00    \\
32 & 3.1250e-02 & 5.6569 & 6.7882e+01  & 8.5513  & 5.0006e-01  \\
64 & 1.5625e-02 & 8.0000  & 9.6000e+01 &8.5184   & 5.0002e-01 \\
128 &7.8125e-03 & 1.1314e+01  & 1.3576e+02  & 8.5018 & 5.0000e-01     \\
%1024 & 9.7656e-04  &   &     \\
\hline
\end{tabular}
\label{sliver2}
\end{table}

\begin{table}[ht]
\caption{$h_{T_0}^3/{|T_0|}$, $H_{T_0}/h_{T_0}$, and $R_3/h_{T_0}$ ($\varepsilon_1 = 1.5$, $\varepsilon_2 = 1.5$)}
\centering
\begin{tabular}{l | l | l | l | l | l } \hline
$N$ &  $s$ & $L_6 / L_1$ & $h_{T_0}^3/{|T_0|}$  & $H_{T_0}/h_{T_0}$ & $R_3/h_{T_0}$\\ \hline \hline
%64 & 1.5625e-02 & 96 & 1.2094e+01   &  1.0000e+00    \\
32 & 3.1250e-02 & 5.6569 & 3.8400e+02 & 3.4986e+01 &  1.4170  \\
64 & 1.5625e-02  & 8.0000  & 7.6800e+02 & 4.8744e+01 & 2.0010 \\
128 &7.8125e-03& 1.1314e+01 & 1.5360e+03 & 6.8411e+01& 2.8288 \\
%1024 & 9.7656e-04  &   &     \\
\hline
\end{tabular}
\label{sliver3}
\end{table}

In Table  \ref{sliver1}, the angle between $\triangle x_1 x_2 x_3$ and $\triangle x_1 x_2 x_4$ tends to $\pi$ as $s \to 0$, and the simplex $T_0$ is ``not good." In Table \ref{sliver2}, the angle between $\triangle x_1 x_3 x_4$ and $\triangle x_2 x_3 x_4$ tends to $0$ as $s \to 0$, and the simplex $T_0$ is ``good." In Table \ref{sliver3}, from the numerical results, the simplex $T_0$ is ``not good."

\begin{figure}[htbp]
\vspace{-4cm}
  \includegraphics[bb=0 0 774 530,scale=0.55]{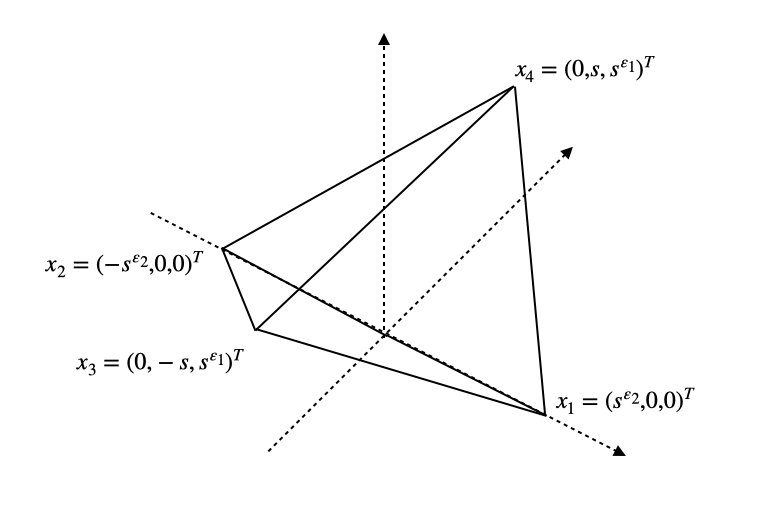}
\caption{Example 8: Sliver}
\label{}
\end{figure}

%We have $\alpha_1 = 2s$, and $\alpha_2 = \alpha_3 = \sqrt{2 s^2 + s^{2 \varepsilon_1}}$, that is,
%\begin{align*}
%\displaystyle
%\frac{\alpha_{\max}}{\alpha_{\min}} \leq c.
%\end{align*}
%On the other hand, in Table \ref{sliver2}, the angle between $\triangle x_1 x_3 x_4$ and $\triangle x_2 x_3 x_4$ tends to $0$ as $s \to 0$. We have $\alpha_1 = \sqrt{2 s^2 + s^{2 \varepsilon_2}}, \alpha_2 = 2 s^{\varepsilon_2}, \alpha_3 = \sqrt{2 s^2 + s^{2 \varepsilon_2}}$, that is,
%\begin{align*}
%\displaystyle
%\frac{\alpha_{\max}}{\alpha_{\min}} \leq c s^{1 - \varepsilon_2} \to \infty \quad \text{as $s \to 0$}.
%\end{align*}
%In the former case, the quantity $\frac{H_T}{h_T}$ diverges, while in the later, the quantity bounds for the constant. These are examples where one can check the validity of the new parameter $H_T$.
\end{Ex}

\begin{comment}
\begin{Ex}[Wedge, Spade, Cap]
In \cite{CheShiZha04}, three types called the wedge, spade, cap, are introduced in addition to the six tetrahedra. The same argument as examples described above is possible for each.
\end{Ex}
\end{comment}

\subsection{Effect of the quantity $|T_0|^{\frac{1}{q} - \frac{1}{p}}$ on the interpolation error estimates for $d = 2,3$} \label{qpinflu}
We now consider the effect of the factor $|T_0|^{\frac{1}{q} - \frac{1}{p}}$. 

\subsubsection{Case in which $q \> p$}
When $q \> p$, the factor may affect the convergence order. In particular, the interpolation error estimate may diverge on anisotropic mesh partitions. 

Let $T_0 \subset \mathbb{R}^2$ be the triangle with vertices $x_1 := (0,0)^T$, $x_2 := (s,0)^T$, $x_3 := (0, s^{\varepsilon})^T$ for $0 \< s \ll 1$, $\varepsilon \geq 1 $, $s \in \mathbb{R}$, and $\varepsilon \in \mathbb{R}$; see Figure \ref{RAT}. Then,
\begin{align*}
\displaystyle
\frac{\alpha_{\max}}{\alpha_{\min}} = s^{1 - \varepsilon}, \quad |T_0| = \frac{1}{2} s^{1 + \varepsilon}.
\end{align*}
%If $\varepsilon = 1$, we use Theorem A, and if $\varepsilon \> 1$, we use Theorem B. 
Let $k=1$, $\ell = 2$,  $m = 1$, $q = 2$, and $ p \in (1,2)$. Then, $W^{1,p}({T}_0) \hookrightarrow L^2({T}_0)$ and Theorem B lead to
\begin{align*}
\displaystyle
| {\varphi}_0 - I_{{T}_0} {\varphi}_0 |_{H^{1}({T}_0)}
\leq  c s^{ - (1+\varepsilon)\frac{2 - p}{2p}} \left( s \left | \frac{\partial \varphi_0}{\partial x_1} \right |_{W^{1,p}({T}_0)} + s^{\varepsilon} \left | \frac{\partial \varphi_0}{\partial x_2} \right |_{W^{1,p}({T}_0)} \right).
\end{align*}
When $\varepsilon = 1$ (i.e., an isotropic element), we obtain
\begin{align*}
\displaystyle
| {\varphi}_0 - I_{{T}_0} {\varphi}_0 |_{H^{1}({T}_0)}
\leq  c h_{T_0}^{ \frac{2(p-1)}{p}} |\varphi_0|_{W^{2,p}(T_0)}, \quad  \frac{2(p-1)}{p} \> 0.
\end{align*}
However, when  $\varepsilon \> 1$ (i.e., an anisotropic element),  the estimate may diverge as $s \to 0$. Therefore, if $q \> p$, the convergence order of the interpolation operator may deteriorate.

\subsubsection{Case in which $q \< p$}
We consider Theorem B. Let $I_{T_0}^{L}: \mathcal{C}^0(T_0) \to \mathcal{P}^k$ ($k \in \mathbb{N}$) be the local Lagrange interpolation operator. Let $  {\varphi}_0  \in W^{\ell,\infty}(T_0) $ be such that $\ell \in \mathbb{N}$, $2 \leq \ell \leq k+1$. Then, for any $m \in \{ 0, \ldots,\ell -1\}$ and $q \in [1,\infty]$, it holds that
\begin{align}
\displaystyle
| {\varphi}_0 - I_{{T}_0}^L {\varphi}_0 |_{W^{m,q}({T}_0)}
\leq c |T_0|^{\frac{1}{q} } \left( \frac{H_{T_0}}{h_{T_0}} \right)^m \sum_{|\gamma| = \ell-m} \mathscr{H}^{\gamma}| \partial^{\gamma}  {(\varphi_0 \circ \Phi_{T_0})} |_{W^{ m ,\infty}(\Phi_{T_0}^{-1}(T_0))}. \label{rem8=56}
\end{align}
Therefore, the convergence order is improved by $|T_0|^{\frac{1}{q} }$.

We can perform some numerical tests to confirm this. Let $k=1$ and 
\begin{align*}
\displaystyle
\varphi_0(x,y,z) := x^2 + \frac{1}{4} y^2 + z^2.
\end{align*}

\begin{comment}
Let $s := \frac{1}{N}$, $N \in \mathbb{N}$ and $\varepsilon \in \mathbb{R}$, $1 \< \varepsilon$. We compute the convergence order with respect to the $H^1$ norm defined by
\begin{align*}
\displaystyle
%&Err_s(H^1) := \frac{| \varphi - I_{T}^{L} \varphi |_{H^1(T)}}{| \varphi |_{H^1(T)}}, 
Err_s^{\varepsilon}(H^1) := | \varphi_0 - I_{T_0}^{L} \varphi_0 |_{H^1(T_0)}.
\end{align*}
The convergence indicator $r$ is defined by
\begin{align*}
\displaystyle
r = \frac{1}{\log(2)} \log \left( \frac{Err_s^{\varepsilon}(H^1)}{Err_{s/2}^{\varepsilon}(H^1)} \right).
\end{align*}
\end{comment}

\begin{description}
  \item[(\Roman{lone})] Let $T_0 \subset \mathbb{R}^3$ be the simplex with vertices $x_1 := (0,0,0)^T$, $x_2 := (s,0,0)^T$, $x_3 := (0,s^{\varepsilon},0)^T$, and $x_4 := (0,0,s^{\delta})^T$ ($1 \< \delta \leq  \varepsilon $), and $0 \< s \ll 1$, $s \in \mathbb{R}$. Then, we have that $\alpha_1 = \sqrt{s^2 + s^{2 \varepsilon}}$, $\alpha_2 = s^{\varepsilon}$, and $\alpha_3 := \sqrt{s^{2 \varepsilon} + s^{2 \delta}}$; i.e.,
\begin{align*}
\displaystyle
\frac{\alpha_{\max}}{\alpha_{\min}}  \leq c s^{1 - \varepsilon}, \quad \frac{H_{T_0}}{h_{T_0}} \leq c.
\end{align*}
From Eq.~\eqref{rem8=56} with $m=1$, $\ell = 2$, and $q=2$, because $|T_0| \approx s^{1+\varepsilon + \delta}$, we have the estimate
\begin{align*}
\displaystyle
&| {\varphi}_0 - I_{{T}_0}^L {\varphi}_0 |_{H^{1}({T}_0)} \leq c  h_{T_0}^{\frac{3 + \varepsilon + \delta}{2}}.
\end{align*}
Computational results for $\varepsilon = 3.0$ and $\delta = 2.0$ are presented in Table \ref{modify=table1}.
  
\begin{table}[htbp]
\caption{Error of the local interpolation operator ($\varepsilon = 3.0, \delta = 2.0$)}
\centering
\begin{tabular}{l | l | l | l  } \hline
$N$ &  $s$  & $Err_s^{3.0}(H^1)$ & $r$   \\ \hline \hline
64 & 1.5625e-02 &  2.4336e-08  &           \\
128 & 7.8125e-03 & 1.5209e-09 &    4.00   \\
256 & 3.9062e-03  & 9.5053e-11  &    4.00      \\
\hline
\end{tabular}
\label{modify=table1}
\end{table}

\begin{figure}[htbp]
\vspace{-5cm}
 \includegraphics[bb=0 0 602 468,scale=0.55]{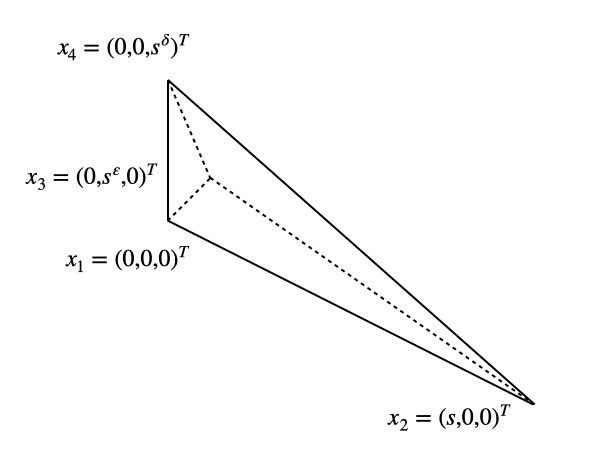}
\caption{Case in which $q \< p$, Example (\Roman{lone})}
\label{}
\end{figure}

\item[(\Roman{ltwo})] Let $T_0 \subset \mathbb{R}^3$ be the simplex with vertices $x_1 := (0,0,0)^T$, $x_2 := (s,0,0)^T$, $x_3 := (s/2,s^{\varepsilon},0)^T$, and $x_4 := (0,0,s)^T$ ($1 \< \varepsilon \leq 6$) and $0 \< s \ll 1$, $s \in \mathbb{R}$. Then, we have that $\alpha_1 = s$, $\alpha_2 = \sqrt{ s^2/4 + s^{2 \varepsilon}}$, and $\alpha_3 := s$; i.e.,
\begin{align*}
\displaystyle
\frac{\alpha_{\max}}{\alpha_{\min}} = \frac{s}{\sqrt{ s^2/4 + s^{2 \varepsilon}}} \leq c, \quad \frac{H_{T_0}}{h_{T_0}} \leq c s^{1 - \varepsilon}.
\end{align*}
From Eq.~\eqref{rem8=56} with $m=1$, $\ell = 2$, and $q=2$, because $|T_0| \approx s^{2+\varepsilon}$, we have the estimate
\begin{align*}
\displaystyle
&| {\varphi}_0 - I_{{T}_0}^L {\varphi}_0 |_{H^{1}({T}_0)} \leq c  h_{T_0}^{3 - \frac{\varepsilon}{2}}.
\end{align*}
Computational results for $\varepsilon = 3.0,6.0$ are presented in Table \ref{modify=table2}.

\begin{table}[htbp]
\caption{Error of the local interpolation operator ($\varepsilon = 3.0,6.0$)}
\centering
\begin{tabular}{l | l | l | l| l| l  } \hline
$N$ &  $s$  & $Err_s^{3.0}(H^1)$ & $r$  & $Err_s^{6.0}(H^1)$ & $r$ \\ \hline \hline
64 & 1.5625e-02 &  1.9934e-04 &      &  1.0206e-01  &     \\
128 & 7.8125e-03 & 7.0477e-05   & 1.50  &  1.0206e-01 & 0\\
256 & 3.9062e-03  & 2.4917e-05  &   1.50 & 1.0206e-01  &  0    \\
\hline
\end{tabular}
\label{modify=table2}
\end{table}
\end{description}

\begin{figure}[htbp]
\vspace{-4cm}
 \includegraphics[bb=0 0 602 468,scale=0.55]{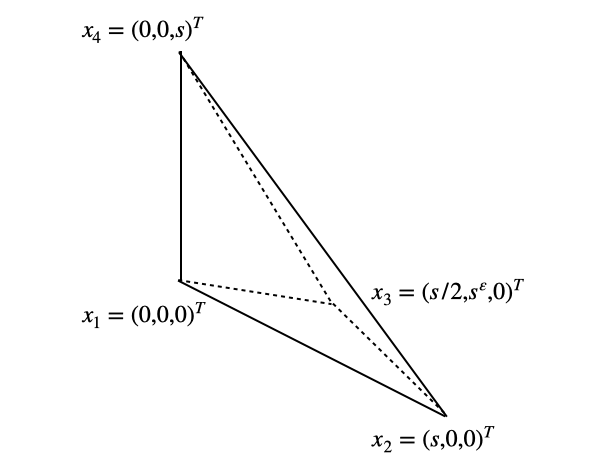}
\caption{Case in which $q \< p$, Example (\Roman{ltwo})}
\label{}
\end{figure}

\subsection{Inverse inequalities} \label{sec_inv}
This section presents some limited results for the inverse inequalities. The results are only stated; the proofs can be found in \cite{Ish22}.

%By analogous argument with the scaling argument, one can easily prove them. Therefore, the results are only stated.

\begin{lem}
Let $\widehat{P} := \mathcal{P}^k$ with $k \in \mathbb{N}$. If Assumption \ref{ass1} is imposed, there exist positive constants $C_i^{IV,d}$, $i=1, \ldots,d$, independent of $h_{T}$ and ${T}$, such that, for all ${\varphi}_h \in{P} = \{  \hat{\varphi}_h \circ {\Phi}_{T}^{-1} ; \ \hat{\varphi}_h \in \widehat{P} \}$,
\begin{align}
\displaystyle
\left \| \frac{\partial \varphi_h}{\partial x_i} \right \|_{L^q(T)}
&\leq C_i^{IV,d} |T|^{\frac{1}{q} - \frac{1}{p}} \frac{1}{\mathscr{H}_i} \| \varphi_h \|_{L^p(T)}. \quad i=1,\ldots,d. \label{inv55}
\end{align}
\end{lem}

\begin{Rem}
If Assumption \ref{ass1} is not imposed, estimate \eqref{inv55} for $i=3$ is
\begin{align}
\displaystyle
\left \| \frac{\partial \varphi_h}{\partial x_3} \right \|_{L^q(T)}
&\leq C_{3}^{IV,3} |T|^{\frac{1}{q} - \frac{1}{p}} \frac{H_{T}}{h_{T}} \frac{1}{\mathscr{H}_2} \| \varphi_h \|_{L^p(T)}.  \label{inv_ex}
\end{align}
This may not be sharp. We leave further arguments for future work.
\end{Rem}

\begin{thrin*}
Let $\widehat{P} := \mathcal{P}^k$ with $k \in \mathbb{N}_0$. Let $\gamma = (\gamma_1,\ldots,\gamma_d) \in \mathbb{N}_0^d$ be a multi-index such that $0 \leq |\gamma| \leq k$. If Assumption \ref{ass1} is imposed, there exists a positive constant $C^{IVC}$, independent of $h_{T}$ and ${T}$, such that, for all ${\varphi}_h \in{P} = \{  \hat{\varphi}_h \circ {\Phi}_{T}^{-1} ; \ \hat{\varphi}_h \in \widehat{P} \}$,
\begin{align}
\displaystyle
\| \partial^{\gamma} \varphi_h \|_{L^q(T)} \leq C^{IVC}  |T|^{\frac{1}{q} - \frac{1}{p}} \mathscr{H}^{- \gamma} \| {\varphi}_h \|_{L^p({T})}. \label{newinv57}
\end{align}
\end{thrin*}

\begin{acknowledgements}
We would like to thank the anonymous referees and the editor of the journal for the valuable comments.
%We thank Stuart Jenkinson, PhD, from Edanz Group (https://jp.edanz.com/ac) for editing a draft of this manuscript. ???
%I have added this acknowledgment to comply with international publishing guidelines (e.g., ICMJE, COPE, EASE) on declaring support given to authors. Please contact us if your target journal asks for a signed letter granting my permission to be acknowledged.

%We thank Professor Norikazu Saito (Tokyo University, Japan) for the useful 	comments.%We are thankful for the time and energy you expended.
%We thank Glenn Pennycook, MSc, from Edanz Group (www.edanzediting.com/ac) for editing a draft of this manuscript.
%If you'd like to thank anyone, place your comments here
%and remove the percent signs.
\end{acknowledgements}

% Authors must disclose all relationships or interests that 
% could have direct or potential influence or impart bias on 
% the work: 
%
% \section*{Conflict of interest}
%
% The authors declare that they have no conflict of interest.

% BibTeX users please use one of
%\bibliographystyle{spbasic}      % basic style, author-year citations
%\bibliographystyle{spmpsci}      % mathematics and physical sciences
%\bibliographystyle{spphys}       % APS-like style for physics
%\bibliography{}   % name your BibTeX data base

% Non-BibTeX users please use

%\newpage
\appendix

\end{document}